\documentclass[12pt,reqno]{amsart}
\usepackage{latexsym}
\usepackage{amssymb}
\usepackage{amscd}

\numberwithin{equation}{section}

\catcode`\@=11
\@namedef{subjclassname@2010}{%
  \textup{2010} Mathematics Subject Classification}
\catcode`\@=12

\catcode`\@=11
\def\FN@{\futurelet\next}
\def\newcodes@{\catcode`\\12\catcode`\{12\catcode`\}12\catcode`\#12%
 \catcode`\%12\relax}
\def\oldcodes@{\catcode`\\0\catcode`\{1\catcode`\}2\catcode`\#6%
 \catcode`\%14\relax}
\def\comment{\newcodes@\endlinechar=10 \comment@}
{\lccode`\0=`\\
\lowercase{\gdef\comment@#1^^J{\comment@@#10endcomment\comment@@@}%
\gdef\comment@@#10endcomment{\FN@\comment@@@}%
\gdef\comment@@@#1\comment@@@{\ifx\next\comment@@@\let\next\comment@
 \else\def\next{\oldcodes@\endlinechar=`\^^M\relax}%
 \fi\next}}}
\catcode`\@=12

\newtheorem{theorem}{Theorem}
\newtheorem*{theoremno}{Theorem}

\newtheorem{lemma}[theorem]{Lemma}

\theoremstyle{remark}
\newtheorem{remark}[theorem]{Remark}
\newtheorem{remarks}[theorem]{Remarks}

\newtheorem*{note}{Note}

\def\al{\alpha}
\def\be{\beta}

\def\la{\lambda}

\def\ph{\varphi}

\def\Z{{\mathbb Z}}

\def\po#1#2{(#1)_#2}

\def\Rat{\operatorname{Rat}}

\def\Pol{\operatorname{Pol}}
\def\fl#1{\left\lfloor#1\right\rfloor}

\def\PSL{\operatorname{PSL}}

\begin{document}
\title[A Riccati differential equation and free subgroup numbers]{A 
Riccati differential equation and free subgroup numbers
for lifts of $\PSL_2(\Z)$ modulo prime powers}
\author[C. Krattenthaler, and 
T.\,W. M\"uller]{C. Krattenthaler$^{\dagger}$ and
T. W. M\"uller$^*$} 

\address{$^{\dagger*}$Fakult\"at f\"ur Mathematik, Universit\"at Wien,
Nordbergstrasze~15, A-1090 Vienna, Austria.
WWW: {\tt http://www.mat.univie.ac.at/\lower0.5ex\hbox{\~{}}kratt}.}

\address{$^*$School of Mathematical Sciences, Queen Mary
\& Westfield College, University of London,
Mile End Road, London E1 4NS, United Kingdom.
}

\thanks{$^\dagger$Research partially supported by the Austrian
Science Foundation FWF, grants Z130-N13 and S9607-N13,
the latter in the framework of the National Research Network
``Analytic Combinatorics and Probabilistic Number
Theory"\newline\indent
$^*$Research supported by Lise Meitner Grant M1201-N13 of the Austrian
Science Foundation FWF}



\begin{abstract}
It is shown that the number $f_\la$ of free subgroups 
of index $6\la$ in the modular group $\PSL_2(\Z)$, when considered
modulo a prime power $p^\al$ with $p\ge5$, is always (ultimately) periodic. 
In fact, an analogous result is established for
a one-parameter family of lifts of
the modular group (containing $\PSL_2(\Z)$ as a special case),
and for a one-parameter family of lifts of the Hecke group
$\mathfrak{H}(4)=C_2*C_4$.
All this is achieved by explicitly determining
Pad\'e approximants to solutions
of a certain multi-parameter family of Riccati differential equations. 
Our main results complement previous work by Kauers and
the authors
([{\it Electron.\ J. Combin.} {\bf 18}(2) (2012), Article~P37]
and [{\it``A
method for determining the mod-$3^k$ behaviour of
recursive sequences"}, preprint]), 
where it is shown, among other things, that the 
free subgroup numbers of $\PSL_2(\Z)$ and 
its lifts display rather complex behaviour modulo powers of $2$ and $3$. 
\end{abstract}
\maketitle

\section{Introduction}
\noindent Beginning with the work of Sylow \cite{Sylow}, Frobenius
\cite{Frob1}, \cite{Frob2}, and P. Hall \cite{PHall1}, \cite{PHall2},
the study of congruences for subgroup numbers and related invariants
has played an important r\^ole in group theory. 
Over the last fifteen years, there has
been active interest in  divisibility properties of subgroup numbers
of (finitely generated) infinite groups; these may, to some extent, be
seen as analogues of classical results for finite groups. The
systematic study of the subgroup arithmetic of infinite groups begins
with \cite{MuHecke}, which investigates the parity of subgroup numbers
and free subgroup numbers in arbitrary Hecke groups $\mathfrak{H}(q) =
C_2\ast C_q$ with $q\geq3$.\footnote{This paper circulated in the 
community since 1998 before finally appearing in print in 2004.} 
The results of \cite{MuHecke} were
subsequently generalised to larger classes of groups and arbitrary
prime modules in \cite{CaMuDescent}, \cite{KratMuHecke},
\cite{MuDescent}, \cite{MuSubArith}, \cite{MuPu2}, and \cite{MuPu}.  

A first attempt at obtaining congruences modulo higher prime powers
was made in \cite{MuPu}, which studies the subgroup numbers of the
inhomogeneous  modular group $\mathfrak{H}(3)\cong \PSL_2(\Z)$
modulo $8$, and derives a congruence modulo $16$ for the number of
free subgroups in $\PSL_2(\Z)$ of given finite index. The real
breakthrough concerning this kind of problem however occurred in
the recent papers \cite{KKM} and \cite{KrMu3Power}, which introduce a
method leading to  semi-automatic existence proofs as well as explicit
computation of such congruences modulo $2$-powers and $3$-powers for a
wide variety of combinatorial sequences. Among other things, this
novel approach leads to congruences modulo arbitrary $2$-powers for
the number $f_\lambda(m)$ of free subgroups of index $6m\lambda$ in
lifts $\Gamma_m(3)$ of $\PSL_2(\Z)$ of the form 
\begin{equation} \label{eq:Gam}
\Gamma_m(3) = C_{2m} \underset{C_m}{\ast} C_{3m} = \big\langle x,
y\,\big\vert\, x^{2m} = y^{3m} = 1,\, x^2=y^3\big\rangle,\quad m\geq1; 
\end{equation}
see Theorems~19 and 20 in \cite{KKM}. The method developed in
\cite{KKM} can in principle be adapted to yield congruences modulo
other prime powers $p^\alpha$, and this has been carried out in detail
in \cite{KrMu3Power} for $p=3$. Among many other results, we obtain
there congruences for the function $f_\lambda(m)$ modulo arbitrary
$3$-powers; cf.\ Section~16 in
\cite{KrMu3Power}. The congruences obtained in this way  
for $f_\lambda(m)$ show a highly non-trivial behaviour of these sequences
modulo $2$-powers and $3$-powers.  
For instance, for the sequence $f_\lambda=f_\lambda(1)$ of 
free subgroup numbers of the group $\mathrm{PSL}_2(\Z)$, one can show that:
\begin{enumerate} 
\item $f_\la\equiv-1~(\text{mod }3)$ if, and only if, 
the $3$-adic expansion of $\la$ is an element of
$\{0,2\}^*1$;
\item $f_\la\equiv1~(\text{mod }3)$ if, and only if, 
the $3$-adic expansion of $\la$ is an element of
$$
\{0,2\}^*100^*\cup \{0,2\}^*122^*;
$$
\item for all other $\la$, we have $f_\la\equiv0~(\text{mod }3)$;
\end{enumerate}
cf.\ \cite[Cor.~53]{KrMu3Power}. Here, for a set $\Omega$, we denote by $\Omega^\ast$ the free monoid generated by $\Omega$.

The general theory of finitely generated virtually free groups
implies that the generating function for the 
free subgroup numbers $f_\la(m)$ satisfies a Riccati differential
equation, see \eqref{eq:fmdiff}. This theory
might also lead one to suspect that the primes $2$ and $3$ should be
special in this context. One of the main results of the present paper
confirms this conjecture in a strong
and surprising form. Namely, we show the following.  

\begin{theoremno}
Let $m$ be a positive integer,
$p$ a prime number with $p\ge 5,$ and let $\al$ be a positive integer. 
Then the generating function $F_m(z)=1+\sum_{\la\ge1}f_\la(m)z^\la,$
when coefficients are reduced modulo $p^\al,$ can be represented as a
rational function. Equivalently, the sequence $(f_\la(m))_{\la\ge1}$
is ultimately periodic modulo $p^\alpha$.
\end{theoremno}

This is Theorem~\ref{thm:Freem} in Section~\ref{Sec:Freem}. The proof
relies on finding an explicit approximate rational solution for a
relevant multi-parameter class of Riccati differential equations; see
Theorem~\ref{thm:ABCD} in the next section. 
In fact, this theorem implies a formula in terms of an explicit
rational function for the generating function $F_m(z)$ when
coefficients are reduced modulo~$p$; see the proof of
Theorem~\ref{thm:Freem} in Section~\ref{Sec:Freem}. A more detailed
$p$-adic analysis of the expressions extracted from
Theorem~\ref{thm:ABCD} allows us to derive even more precise
results on the form of the rational functions and on the period
lengths of the sequence $(f_\la(m))_{\la\ge1}$ modulo concrete prime
powers $p^\alpha$; see Section~\ref{Sec:Periods}. Moreover,
Theorem~\ref{thm:ABCD} is general enough to also produce analogous
results for the free subgroup numbers of lifts of the Hecke group
$\mathfrak{H}(4)$; see Section~\ref{Sec:Freem4}.

The proof of
Theorem~\ref{thm:ABCD} occupies
Section~\ref{Sec:ABCDProof}. Section~\ref{Sec:Freem} presents the
(very short)  derivation of the theorem above from
Theorem~\ref{thm:ABCD}. As already said, 
Section~\ref{Sec:Periods} establishes
more precise results for the subgroup numbers $f_\lambda(m)$.
More specifically, Theorem~\ref{thm:Freep} identifies the denominators of the 
rational functions which one obtains when
the coefficients of the generating function
$1+\sum_{\la\ge1}f_\la(m) z^\la$ are reduced modulo a prime power
$p^\al$ with $p\ge5$. Theorems~\ref{thm:Free7}--\ref{thm:Free13^5}
then illustrate this general result for the numbers 
$f_\lambda=f_\lambda(1)$ (i.e.,
the free subgroup numbers of the modular group $\PSL_2(\Z)$ itself) 
when reduced modulo the primes
$p=7,11,13$. In particular, the exact minimal period is determined
in these cases. Section~\ref{Sec:Freem4} presents the analogous
results for the free subgroup numbers of the
lifts of $\mathfrak{H}(4)$.
Finally, the rather interesting and involved way which led us
to conjecture formulae (\ref{eq:pnk}) and (\ref{eq:qnk}) for numerator
and denominator of the approximate rational solution $R_n(z)$ of the
differential equation (\ref{eq:diffeqABCD}) 
(which is the subject of Theorem~\ref{thm:ABCD})
is explained in an appendix. 
 
\begin{note}
This paper is accompanied by two {\sl Mathematica} files
that allow an interested reader to compute the rational functions
representing the reductions modulo $p^\alpha$ of the generating 
function $1+\sum_{\la\ge1}f_\la z^\la$
for the number $f_\la$ of free subgroups
of index $6\la$ in $\PSL_2(\Z)$ and of the generating 
function $1+\sum_{\la\ge1}f_\la^{(4)} z^\la$
for the number $f_\la^{(4)}$ of free subgroups
of index $4\la$ in $\mathfrak{H}(4)$.
The result is given as partial fraction expansion of the form
predicted by Theorems~\ref{thm:Freep} and \ref{thm:Freep4},
respectively.
The files are available at the article's
website
{\tt http://www.mat.univie.ac.at/\lower0.5ex\hbox{\~{}}kratt/artikel/psl2zmod.html}.
\end{note}

\section{Approximate rational solutions for a Riccati differential equation}

\noindent We consider the Riccati differential equation
\begin{equation} \label{eq:diffeqABCD}
(1-Az)F(z)-Bz^2F'(z)-CzF^2(z)-1-Dz=0,
\end{equation}
where $A,B,C,D$ are constants.
The technical main result of our paper is the following. 

\begin{theorem} \label{thm:ABCD}
For every positive integer $n,$
there exist uniquely determined 
polynomials $P_n(z)=1+\sum_{k=1}^n p_{n,k}z^k$ and 
$Q_n(z)=1+\sum_{k=1}^n q_{n,k}z^k,$ where $p_{n,k}$ and $q_{n,k}$
are homogeneous polynomials in $A,B,C,D$ of degree $k$
over the integers, such that 
the rational function $R_n(z)=P_n(z)/Q_n(z)$ satisfies the
differential equation 
\begin{multline} \label{eq:idABCD}
(1-Az)R_n(z)-Bz^2R_n'(z)-CzR_n^2(z)-1-Dz\\
=
-\frac {z^{2n+1}} {Q_n^2(z)}(A+C+D)
\prod _{\ell=1} ^{n}(\ell AB+AC+CD+\ell^2B^2+2\ell BC+C^2).
\end{multline}
Moreover, the coefficients $p_{n,k}$ and $q_{n,k}$ are given explicitly by
\begin{multline} \label{eq:pnk}
p_{n,n-k}=\frac {(-1)^{n+1}B^{n-k}} {2C\,(\frac {E} {B}-k)_{2k+1}}
\\
\times
\Bigg( \left(\tfrac {A + 2 C + E} {2B}\right)_{n+1}
\sum_{j=0}^k \binom {k+j}k\binom {n-j}{k-j}
\left(-\tfrac E B
  + j + 1\right)_{ 
    k - j} \left(\tfrac {A + 2 C - E} {2B}\right)_j\\
\cdot
    \big (A +\tfrac{ 2 k j}{k+j} B - \tfrac{k - j}{k+j} E\big)\\
- \left(\tfrac {A + 2 C - E} {2B}\right)_{n+1}
\sum_{j=0}^k \binom {k+j}k\binom {n-j}{k-j} 
\left(\tfrac E B
  + j + 1\right)_{ 
    k - j} \left(\tfrac {A + 2 C + E} {2B}\right)_j\\
\cdot
    \big (A +\tfrac{ 2 k j}{k+j} B + \tfrac{k - j}{k+j} E\big)
\Bigg)
\end{multline}
and
\begin{multline} \label{eq:qnk}
q_{n,n-k}=\frac {(-1)^nB^{n-k}} {(\frac {E} {B}-k)_{2k+1}}
\\
\times
\Bigg( \left(\tfrac {A + 2 C + E} {2B}\right)_{n+1}
\sum_{j=0}^k \binom {k+j}k\binom {n-j}{k-j}
\left(-\tfrac E B
  + j + 1\right)_{ 
    k - j} \left(\tfrac {A + 2 C - E} {2B}\right)_j\\
- \left(\tfrac {A + 2 C - E} {2B}\right)_{n+1}
\sum_{j=0}^k \binom {k+j}k\binom {n-j}{k-j}
\left(\tfrac E B
  + j + 1\right)_{ 
    k - j} \left(\tfrac {A + 2 C + E} {2B}\right)_j
\Bigg),
\end{multline}
where $E^2=A^2-4CD,$ and
the {\em Pochhammer symbol} $(\al)_m$ is defined by
$(\alpha)_m  =  \alpha(\hbox{$\alpha+1$})(\alpha+2)\cdots (\alpha+m-1)$ for
$m>0,$ and $(\alpha)_0  = 1$.
\end{theorem}

\begin{remark} \label{rem:1}
(1) In the sums in \eqref{eq:pnk} it might seem that there is a problem if
  $k=j=0$. 
The correct way to interpret \eqref{eq:pnk} for $k=0$ is the
  following: one first leaves $k$ undetermined; then there is no
problem to calculate the expression
$A +\tfrac{ 2 k j}{k+j} B - \tfrac{k - j}{k+j} E$
when $j=0$, it simply equals $A-E$. 
Subsequently, one can safely set $k=0$.
A similar remark applies to its
``companion expression" in the second sum. So,
explicitly,
\begin{equation} \label{eq:pn0} 
p_{n,n}=\frac {(-1)^{n+1}B^{n+1}} {2CE}
\Big(    (A  -  E)
 \left(\tfrac {A + 2 C + E} {2B}\right)_{n+1}
-(A+E) \left(\tfrac {A + 2 C - E} {2B}\right)_{n+1}
\Big).
\end{equation}

\medskip
\noindent (2) The rational function solutions $R_n(z)$ for $n=1,2,3$
are displayed in the appendix as  
\eqref{eq:R1}--\eqref{eq:R3}.

\medskip
\noindent (3) Formulae for the first few coefficients $p_{n,k}$ and $q_{n,k}$,
$k=1,2,3$, are displayed in \eqref{eq:pn1}--\eqref{eq:qn3}
in the appendix.
Assuming Theorem~\ref{thm:ABCD}, these (originally experimentally
found) formulae can be established as follows: according to the theorem, 
the coefficients $p_{n,k}$ and $q_{n,k}$ are in particular polynomials
in $A$ of degree at most $k$. As such, they are uniquely determined
by $k+1$ special evaluations. These can be found from Formulae~\eqref{eq:pnk}
and \eqref{eq:qnk} by setting $A=-2C-E-2Bs$, $s=0,1,\dots,k$.
Indeed, with these specialisations, the term
$\left(\tfrac {A + 2 C + E} {2B}\right)_{n+1}$ vanishes for $n\ge k$,
so that the first terms between parentheses in \eqref{eq:pnk} and 
\eqref{eq:qnk} vanish. For the same reason,
the term $\left(\tfrac {A + 2 C + E} {2B}\right)_{j}$ in the second sums
over $j$ in \eqref{eq:pnk} and \eqref{eq:qnk} vanishes for $j\ge k+1$,
so that only the terms for $j=0,1,\dots,k$ remain.

\medskip
\noindent (4)
It is easy to see that
$R_n(z)$ solves \eqref{eq:idABCD} if, and only if, 
$R_n(z/B)$ solves \eqref{eq:idABCD} with 
$B$ replaced by $1$,
$A$ replaced by $A/B$, 
$C$ replaced by $C/B$, and 
$D$ replaced by $D/B$. 
Hence, it suffices to consider the special case where $B=1$.
We shall occasionally make use of this simplification of the problem.
\medskip

\noindent (5) It is in fact not necessary to specify the precise form of the
right-hand side in \eqref{eq:idABCD} once we know that it is of
the form $\mathcal O(z^{2n+1})$, that is, a formal power series of
order at least $2n+1$; see Remark~\ref{rem:eind}.
Consequently, a different way to view the rational function $R_n(z)$ 
is as the {\it Pad\'e approximant\/} of order $n$ to the solution of the
Riccati differential equation \eqref{eq:diffeqABCD}. See also
the proof of Lemma~\ref{lem:pqeind}.
\end{remark}

\section{Proof of Theorem~\ref{thm:ABCD}}
\label{Sec:ABCDProof}

\noindent The proof of Theorem~\ref{thm:ABCD} is accomplished by
dividing the assertion of the theorem into several parts,  
and then treating each of these partial assertions  
in a separate lemma. More precisely, in
Lemma~\ref{lem:pqsat} below we show that the rational function
$R_n(z)=P_n(z)/Q_n(z)$ given by \eqref{eq:pnk} and \eqref{eq:qnk}
does indeed satisfy the differential equation \eqref{eq:idABCD}, while  
uniqueness of solution is the subject of Lemma~\ref{lem:pqeind}. 
The fact that the coefficients
$p_{n,k}$ and $q_{n,k}$ in \eqref{eq:pnk} and \eqref{eq:qnk} are
homogeneous polynomials in $A,B,C,D$ of degree $k$ is then established
in Lemma~\ref{lem:pqpol}, while 
integrality of their coefficients is proved in Lemma~\ref{lem:pqint}.
Taken together, these 
lemmas provide a full proof of Theorem~\ref{thm:ABCD}.

\begin{lemma} \label{lem:pqsat}
The rational function $R_n(z)=P_n(z)/Q_n(z),$ with the coefficients of
the polynomials $P_n(z)$ and $Q_n(z)$ given by \eqref{eq:pnk} and
\eqref{eq:qnk}{\em ,} satisfies the differential equation \eqref{eq:idABCD}.
\end{lemma}

\begin{proof}
If we substitute $R_n(z)=P_n(z)/Q_n(z)$ into the differential equation
and multiply both sides by $Q_n^2(z)$, then we see that, in order to
establish the lemma, we have to
demonstrate the equation
\begin{multline} \label{eq:idPQ} 
(1-Az)P_n(z)Q_n(z)-Bz^2(P_n'(z)Q_n(z)-P_n(z)Q_n'(z))
-CzP_n^2(z)-(1+Dz)Q_n^2(z)\\
=
-z^{2n+1}(A+C+D)
\prod _{\ell=1} ^{n}(\ell AB+AC+CD+\ell^2B^2+2\ell BC+C^2).
\end{multline}

We start by verifying that the coefficients of $z^{2n+1}$ 
on both sides of \eqref{eq:idPQ} agree. For the sake of better
readability, let us write
\begin{equation} \label{eq:Pi} 
\Pi_+:=B^{n+1}\left(\tfrac {A + 2 C + E} {2B}\right)_{n+1}\quad
\text{and}\quad 
\Pi_-:=B^{n+1}\left(\tfrac {A + 2 C - E} {2B}\right)_{n+1}.
\end{equation}
Using this notation plus \eqref{eq:pn0}, the coefficient
$p_{n,n}$ is seen to be given by
\begin{equation} \label{eq:pnn} 
p_{n,n}=\frac {(-1)^{n+1}} {2CE}
\big(    (A  -  E)
 \Pi_+
-(A+E)\Pi_-
\big),
\end{equation}
while, by \eqref{eq:qnk}, the coefficient $q_{n,n}$ is given by
\begin{equation} \label{eq:qnn} 
q_{n,n}=\frac {(-1)^n} {E}
\big( \Pi_+
- \Pi_-
\big).
\end{equation}
Consequently, the coefficient of $z^{2n+1}$ on the left-hand side of
\eqref{eq:idPQ} equals
\begin{align*}
-Ap_{n,n}q_{n,n}
-Cp^2_{n,n}
-Dq^2_{n,n}&=
\frac {1} {4C^2E^2}
\Big(
\big(A(A-E)(2C)-C(A-E)^2-D(2C)^2\big)\Pi_+^2\\
&\kern1cm
+\big(-2A^2(2C)+2C(A^2-E^2)+2D(2C)^2\big)\Pi_+\Pi_-\\
&\kern1cm
+\big(A(A+E)(2C)-C(A+E)^2-D(2C)^2\big)\Pi_-^2
\Big)\\
&=-\frac {1} {C}\Pi_+\Pi_-,
\end{align*}
where we have used the relation $E^2=A^2-4CD$ in the last step. 
It is straightforward to see that this is identical with the
coefficient of $z^{2n+1}$ on the right-hand side of \eqref{eq:idPQ}.

We now divide both sides of \eqref{eq:idPQ} by $z^{2n+1}$,
then differentiate them with respect to $z$, and finally multiply
both sides of the resulting equation by $z^{2n+2}$. In this way, 
we obtain the equation 
\begin{multline} \label{eq:idPQ2}
P_n(z) \big(2 C n z
   P_n(z)-2 C z^2
   P_n'(z)-(2 n+1) Q_n(z)+z
   Q_n'(z)\\
-z^2 (A+2 B n-B)
   Q_n'(z)+2 A n z Q_n(z)+B z^3
   Q_n''(z)\big)\\
+Q_n(z)
   \big(z P_n'(z)-z^2 (A-2 B n+B)
   P_n'(z)-B z^3 P_n''(z)\kern3cm\\
+(2 n+1)
   Q_n(z)+2
   D n z Q_n(z)-2
   D z^2 Q_n'(z)-2 z
   Q_n'(z)\big)=0.
\end{multline}
Since we already know that the coefficients of $z^{2n+1}$ on both
sides of \eqref{eq:idPQ} agree, establishing \eqref{eq:idPQ2} is
completely equivalent to establishing \eqref{eq:idPQ}. We concentrate
on the former task from now on.

The crucial observation is that, if \eqref{eq:idPQ2} is to hold,
then the first factor between parentheses in \eqref{eq:idPQ2} must be a
(polynomial) multiple of 
$Q_n(z)$, while the second factor between
parentheses must be a multiple of 
$P_n(z)$. Indeed, our claim is that
\begin{multline} \label{eq:Pfaktor} 
2 C n z
   P_n(z)-2 C z^2
   P_n'(z)-(2 n+1) Q_n(z)+z
   Q_n'(z)
-z^2 (A+2 B n-B)
   Q_n'(z)\\
+2 A n z Q_n(z)+B z^3
   Q_n''(z)
= -(2n+1-nAz+n^2Bz)Q_n(z)
\end{multline}
and that 
\begin{multline} \label{eq:Qfaktor} 
z P_n'(z)-z^2 (A-2 B n+B)
   P_n'(z)-B z^3 P_n''(z)
+(2 n+1)
   Q_n(z)+2
   D n z Q_n(z)\\
-2
   D z^2 Q_n'(z)-2 z
   Q_n'(z)
=(2n+1-nAz+n^2Bz) P_n(z).
\end{multline}
Clearly, if we manage to establish \eqref{eq:Pfaktor} and
\eqref{eq:Qfaktor}, then \eqref{eq:idPQ2}, and hence also \eqref{eq:idPQ}
and \eqref{eq:idABCD} follow immediately.

Let us write
\begin{align*}
p_{n,k}&=p_{n,k}^+\Pi_+ -p_{n,k}^-\Pi_-,\\
q_{n,k}&=q_{n,k}^+\Pi_+ -q_{n,k}^-\Pi_-,
\end{align*}
where $\Pi_+$ and $\Pi_-$ are given by \eqref{eq:Pi} as before,
which implicitly defines the $p_{n,k}^+$'s, etc., and let
\begin{alignat*}2
P_n^+(z)&=\sum_{k=0}^n p_{n,k}^+z^k,
&P_n^-(z)&=\sum_{k=0}^n p_{n,k}^-z^k,\\
Q_n^+(z)&=\sum_{k=0}^n q_{n,k}^+z^k,\quad \quad 
&Q_n^-(z)&=\sum_{k=0}^n q_{n,k}^-z^k.
\end{alignat*}
One might hope that \eqref{eq:Pfaktor} and \eqref{eq:Qfaktor} continue
to hold if one replaces $P_n(z)$ by $P_n^+(z)$ and
$Q_n(z)$ by $Q_n^+(z)$ everywhere, and that they continue to hold
if one replaces $P_n(z)$ by $P_n^-(z)$ and
$Q_n(z)$ by $Q_n^-(z)$ everywhere. This turns out to be too
optimistic, but not too far from the truth.
Our new claim is that
\begin{multline} \label{eq:Pfaktor+} 
2 C n z
   P_n^+(z)-2 C z^2
   (P_n^+)'(z)-(2 n+1) Q_n^+(z)+z
   (Q_n^+)'(z)
-z^2 (A+2 B n-B)
   (Q_n^+)'(z)\\
+2 A n z Q_n^+(z)+B z^3
   (Q_n^+)''(z)
= -(2n+1-nAz+n^2Bz)Q_n^+(z),
\end{multline}
\begin{multline} \label{eq:Pfaktor-} 
2 C n z
   P_n^-(z)-2 C z^2
   (P_n^-)'(z)-(2 n+1) Q_n^-(z)+z
   (Q_n^-)'(z)
-z^2 (A+2 B n-B)
   (Q_n^-)'(z)\\
+2 A n z Q_n^-(z)+B z^3
   (Q_n^-)''(z)
= -(2n+1-nAz+n^2Bz)Q_n^-(z),
\end{multline}
\begin{multline} \label{eq:Qfaktor+} 
z (P_n^+)'(z)-z^2 (A-2 B n+B)
   (P_n^+)'(z)-B z^3 (P_n^+)''(z)
+(2 n+1)
   Q_n^+(z)+2
   D n z Q_n^+(z)\\
-2
   D z^2 (Q_n^+)'(z)-2 z
   (Q_n^+)'(z)
=(2n+1-nAz+n^2Bz) P_n^+(z)\\
+(2n+1)\binom {2n}n
\frac {(-1)^nB^{-n}\, 
 \left(\tfrac {A + 2 C - E} {2B}\right)_{n+1}} 
{C\,(\frac {E} {B}-n)_{2n+1}},
\end{multline}
\begin{multline} \label{eq:Qfaktor-} 
z (P_n^-)'(z)-z^2 (A-2 B n+B)
   (P_n^-)'(z)-B z^3 (P_n^-)''(z)
+(2 n+1)
   Q_n^-(z)+2
   D n z Q_n^-(z)\\
-2
   D z^2 (Q_n^-)'(z)-2 z
   (Q_n^-)'(z)
=(2n+1-nAz+n^2Bz) P_n^-(z)\\
+(2n+1)\binom {2n}n
\frac {(-1)^nB^{-n}\, 
 \left(\tfrac {A + 2 C + E} {2B}\right)_{n+1}} 
{C\,(\frac {E} {B}-n)_{2n+1}}.
\end{multline}

These four equations do indeed imply \eqref{eq:Pfaktor} and
\eqref{eq:Qfaktor}: to be 
precise, by multiplying both sides of
\eqref{eq:Pfaktor+} by $\Pi_+$ and multiplying both sides of
\eqref{eq:Pfaktor-} by $\Pi_-$, and then taking the difference of the
two equations, we obtain \eqref{eq:Pfaktor}. An analogous argument
shows that \eqref{eq:Qfaktor+} and \eqref{eq:Qfaktor-} together imply
\eqref{eq:Qfaktor}. Moreover, since
$$P_n^+(z)\big\vert_{E\to-E}=-P_n^-(z)\quad \text{and}\quad 
Q_n^+(z)\big\vert_{E\to-E}=-Q_n^-(z),$$
Equation~\eqref{eq:Pfaktor-} results from \eqref{eq:Pfaktor+} by
replacing $E$ by $-E$ (note that the parameter $D=(A^2-E^2)/(4C)$
remains invariant under this replacement), and  
an analogous assertion
holds for \eqref{eq:Qfaktor+} and \eqref{eq:Qfaktor-}. Thus, the proof
of the lemma has been reduced to establishing  \eqref{eq:Pfaktor+} and
\eqref{eq:Qfaktor+}. 

We start with \eqref{eq:Pfaktor+}. Both sides of this equation are
polynomials in $z$ and $A$ (with coefficients depending on $n,B,C,E$).
Due to the form of the coefficients $p_{n,k}^+$ and $q_{n,k}^+$, every
subexpression in \eqref{eq:Pfaktor+} may conveniently (if 
tediously) be expanded in the basis
$$
\big\{z^k\left(\tfrac {A + 2 C - E} {2B}\right)_j:j,k=0,1,\dots\big\}.
$$
It can be readily verified that corresponding coefficients
on both sides of \eqref{eq:Pfaktor+} agree.

We use the same approach for \eqref{eq:Qfaktor+}. Again, one readily
verifies that corresponding coefficients
on both sides of \eqref{eq:Qfaktor+} agree, except possibly for the
coefficients of basis elements of the
form $z^0\left(\tfrac {A + 2 C - E} {2B}\right)_j$, $j=0,1,\dots$
In other words, we have shown that \eqref{eq:Qfaktor+} holds except
possibly for an error in the coefficient of $z^0$. There are only three
terms in \eqref{eq:Qfaktor+} which have non-zero coefficient of $z^0$:
the term $(2 n+1)Q_n^+(z)$ on the left-hand side, the term
$(2n+1)P_n^+(z)$ on
the right-hand side, and the last term on the right-hand side. 
So, in order to completely establish
\eqref{eq:Qfaktor+}, we have to show that
$$
(2n+1)q_{n,0}^+-(2n+1)p_{n,0}^+=
(2n+1)\binom {2n}n
\frac {(-1)^nB^{-n}\, 
 \left(\tfrac {A + 2 C - E} {2B}\right)_{n+1}} 
{C\,(\frac {E} {B}-n)_{2n+1}},
$$
or, explicitly, that 
\begin{multline} \label{eq:Gosperid}
(2n+1)
\frac {(-1)^nB^{-n-1}} {2C\,(\frac {E} {B}-n)_{2n+1}}
\sum_{j=0}^n \binom {n+j}n
\left(-\tfrac E B
  + j + 1\right)_{ 
    n - j} \left(\tfrac {A + 2 C - E} {2B}\right)_j\\
\cdot
    \big (2C+A +\tfrac{ 2 n j}{n+j} B - \tfrac{n - j}{n+j} E\big)
=
(2n+1)\binom {2n}n
\frac {(-1)^nB^{-n}\, 
 \left(\tfrac {A + 2 C - E} {2B}\right)_{n+1}} 
{C\,(\frac {E} {B}-n)_{2n+1}}.
\end{multline}
The Gosper algorithm (cf.\ 
\cite{GospAB}, \cite[\S~5.7]{GrKPAA}, \cite[\S~II.5]{PeWZAA}) finds that
\begin{multline} \label{eq:Gosper}
(2n+1)
\frac {(-1)^nB^{-n-1}} {2C\,(\frac {E} {B}-n)_{2n+1}}
 \binom {n+j}n
\left(-\tfrac E B
  + j + 1\right)_{ 
    n - j} \left(\tfrac {A + 2 C - E} {2B}\right)_j\\
\cdot
    \big (2C+A +\tfrac{ 2 n j}{n+j} B - \tfrac{n - j}{n+j} E\big)
=G(n,j+1)-G(n,j),
\end{multline}
with
$$G(n,j)=
(2n+1)
\frac {(-1)^nB^{-n}} {C\,(\frac {E} {B}-n)_{2n+1}}
\binom {n+j-1}n
\left(-\tfrac E B
  + j \right)_{ 
    n - j+1} \left(\tfrac {A + 2 C - E} {2B}\right)_j.
$$
Clearly, summing both sides of \eqref{eq:Gosper} over $j$ running
from $0$ to $n$ immediately yields \eqref{eq:Gosperid}, and thus
\eqref{eq:Qfaktor+}.

This completes the proof of the lemma.
\end{proof}

\begin{lemma} \label{lem:pqeind}
The differential equation \eqref{eq:idABCD} has a unique solution 
$R_n(z)=P_n(z)/Q_n(z),$ where $P_n(z)$ and $Q_n(z)$ are polynomials in
$z$ of degree at most $n,$ and $P_n(0)=Q_n(0)=1$.
\end{lemma}

\begin{proof}
The first observation is that the differential equation
\eqref{eq:diffeqABCD} has a unique formal power series solution
$F(z)$. This is seen by comparing
coefficients of $z^m$, $m=0,1,\dots$, which allows one to compute the
coefficients of the series $F(z)$ recursively. More precisely,
comparison of constant coefficients yields that $F(0)=1$, that is,
that the coefficient of $z^0$ in $F(z)$ equals~$1$, while, for $m\geq1$, 
comparison of coefficients of $z^m$ yields an equation
which expresses the coefficient of $z^m$ in $F(z)$ uniquely in terms
of lower (already computed) coefficients.

Next one observes that
the differential equations \eqref{eq:diffeqABCD} and \eqref{eq:idABCD}
agree for the coefficients of $z^m$, $m=0,1,\dots,2n$. Thus, what we are
computing in \eqref{eq:idABCD} is the {\it Pad\'e approximant} 
of the uniquely determined solution $F(z)$ of \eqref{eq:diffeqABCD}
with the degrees of the numerator and denominator 
polynomials bounded above by $n$.
Since Pad\'e approximants are unique (see e.g.\ \cite{BaGrAA}),
a rational function $R_n(z)$ solving \eqref{eq:idABCD}, if there is
one at all, 
is uniquely determined. However, from Lemma~\ref{lem:pqsat}, we know
that such a solution does indeed exist; and, by the arguments above, it is 
the unique solution of \eqref{eq:idABCD}.
\end{proof}

\begin{remark} \label{rem:eind}
It should be observed that the above proof of Lemma~\ref{lem:pqeind}
did not use the precise form of the right-hand side of
\eqref{eq:idABCD}. Rather it shows that, if there is a rational
function solution of \eqref{eq:idABCD} of the prescribed form, 
then this right-hand side is {\it forced\/} upon us.
\end{remark}

\begin{lemma} \label{lem:pqpol}
The coefficients $p_{n,k}$ and $q_{n,k},$ as defined in \eqref{eq:pnk}
and \eqref{eq:qnk}{\em ,} are homogeneous polynomials in $A,B,C,D$ of
degree $k$. 
In particular, $p_{n,0}=q_{n,0}=1$.
\end{lemma}
\begin{proof}
We first consider $q_{n,k}$. We shall prove that the expression in
\eqref{eq:qnk} is a polynomial in $A,B,C,E$. 
Taking into account that the denominator can be rewritten as
\begin{equation} \label{eq:Eprod} 
\left(\tfrac {E} {B}-k\right)_{2k+1}=E
\prod _{i=1} ^{k}\left(\tfrac {E^2} {B^2}-i^2\right),
\end{equation}
and that the term between parentheses in \eqref{eq:qnk} is of the form
$\Pol(E)-\Pol(-E)$, where $\Pol(E)$ is some polynomial in $E$,
we see that $q_{n,k}$ is a rational function in $E^2$. Hence,
if we are able to show that $q_{n,k}$ is a
{\it polynomial\/} in $A,B,C,E$,
then, in view of the relation $E^2=A^2-4CD$, it will be obvious that
$q_{n,k}$ is also a
polynomial in $A,B,C,D$.

In order to accomplish this,
we shall show that the expression between parentheses in
\eqref{eq:qnk} vanishes 
for $E=Bs$ with $s=-k,-k+1,\dots,k$. This proves that this expression
is divisible by the denominator on the right-hand side of \eqref{eq:qnk}
(rewritten here in \eqref{eq:Eprod}) and, after
multiplying numerator and denominator of the expression in
\eqref{eq:qnk} by $B^{2k+1}$, 
that $q_{n,k}$ is a homogeneous
polynomial in $A,B,C,E$ of degree $k$, as desired. 

Writing the first term between parentheses, with $E=Bs$, 
in standard hypergeometric notation, 
$${}_p F_q\!\left[\begin{matrix} a_1,\dots,a_p\\ b_1,\dots,b_q\end{matrix}; 
z\right]=\sum _{m=0} ^{\infty}\frac {\po{a_1}{m}\cdots\po{a_p}{m}}
{m!\,\po{b_1}{m}\cdots\po{b_q}{m}} z^m\ ,
$$
we get 
\begin{multline} \label{eq:term1}
\left(\tfrac {A + 2 C + Bs} {2B}\right)_{n+1}
\sum_{j=0}^k \binom {k+j}k\binom {n-j}{k-j}
\left(-s
  + j + 1\right)_{ 
    k - j} \left(\tfrac {A + 2 C - Bs} {2B}\right)_j\\
=
 \left(\tfrac {A + 2 C + Bs} {2B}\right)_{n+1}\,
   \left(\tfrac {A + 2 C - Bs} {2B}\right)_s\,
\binom {k+s}k\frac {(n-s)!} {(n-k)!}\\
\times
{}_3 F_2\!\left[\begin{matrix}
-k+s,k+s+1,\frac {A} {2B}+\frac {C} {B}+\frac {s} {2}\\
-n+s,s+1
\end{matrix};1\right].
\end{multline}
On the other hand, if we write the second term between parentheses
with $E=Bs$ in hypergeometric notation, then we get
\begin{multline} \label{eq:term2}
 \left(\tfrac {A + 2 C - Bs} {2B}\right)_{n+1}
\sum_{j=0}^k \binom {k+j}k\binom {n-j}{k-j}
\left(s
  + j + 1\right)_{ 
    k - j} \left(\tfrac {A + 2 C + Bs} {2B}\right)_j\\
= 
\left(\tfrac {A + 2 C - Bs} {2B}\right)_{n+1}\,
\binom nk
\frac {(k+s)!} {s!}
\,{}_3 F_2\!\left[\begin{matrix}
-k,k+1,\frac {A} {2B}+\frac {C} {B}+\frac {s} {2}\\
-n,s+1
\end{matrix};1\right].
\end{multline}
In order to see that one expression can be transformed into the
other, we apply the transformation formula (see \cite[Ex.~7, p.~98]{BailAA})
\begin{equation}
{} _{3} F _{2} \!\left [ \begin{matrix} {a, b, c}
\\ {d, e}\end{matrix};1
\right ] =\frac {\Gamma(e)\,\Gamma(d+e-a-b-c)}
{\Gamma(e-a)\,\Gamma( d+e-b-c)}
{} _{3} F _{2} \!\left [ \begin{matrix} {a, d-b, d-c}
\\ {d, d+e-b-c}\end{matrix};1
\right ]
\label{eq:T3204}
\end{equation}
with $a=\frac {A} {2B}+\frac {C} {B}+\frac {s} {2}$, 
$b=-k$, $c=k+1$, $d=s+1$, $e\to-n$ to \eqref{eq:term2} 
(here, $\Gamma(x)$ denotes the classical gamma function). After some
manipulation, one sees that the result is exactly \eqref{eq:term1}.

\medskip
The arguments for $p_{n,k}$ are completely analogous, albeit
more involved due to existence of the extra linear terms in 
$A,C,E$ in the sums.
One additional detail here is that one also has to show that $C$
divides the expression between parentheses in \eqref{eq:pnk}.
In order to see this, one sets $C=0$ in this expression.
In view of $E^2=A^2-4CD$, this specialisation implies that
$E=\pm A$. Without loss of generality, we let $E=A$. Then 
$ \left(\tfrac {A + 2 C - E} {2B}\right)_j=0$ for all $j>0$.
Hence, the only term possibly surviving is the one for $j=0$ in the
first sum between 
parentheses in \eqref{eq:pnk}. 
However, for $j=0$, $C=0$, and $E=A$ we have 
$A +\tfrac{ 2 k j}{k+j} B - \tfrac{k - j}{k+j} E=
A  - \tfrac{k }{k} A=0$, so that this term vanishes also.

Now, that we already know that $p_{n,0}$ and $q_{n,0}$ are constants,
that is, that they are 
independent of $A,B,C,D$ (or, equivalently, of $A,B,C,E$),
in order to see that in fact $p_{n,0}=q_{n,0}=1$, we are free to
arbitrarily specialise the variables $A,B,C,E$. Our choice is
$A=-2C-n-1$, $B=1$, and $E=-n-1$. Then we have
$ \left(\tfrac {A + 2 C - E} {2B}\right)_{n+1}=
 (0)_{n+1}=0$, and therefore, in both \eqref{eq:pnk} and \eqref{eq:qnk},
the second term between parentheses vanishes. For the same reason,
in the remaining sums over $j$ only the terms for $j=0$ 
survive, so that we obtain 
$$
p_{n,0}=\frac {(-1)^{n+1}} {2C\,(-2n-1)_{2n+1}}
\left(-n-1 \right)_{n+1}
\left(n 
  +  2\right)_{ 
    n }
     (-2C)=1
$$
and, similarly, 
$$
q_{n,0}=\frac {(-1)^{n}} {(-2n-1)_{2n+1}}
\left(-n-1 \right)_{n+1}
\left(n 
  +  2\right)_{ 
    n }
    =1.
$$

This completes the proof of the lemma.
\end{proof}

\begin{lemma} \label{lem:pqint}
The coefficients $p_{n,k}$ and $q_{n,k},$ as defined in \eqref{eq:pnk}
and \eqref{eq:qnk}{\em ,} are polynomials in $A,B,C,D$ with integer
coefficients. 
\end{lemma}

\begin{proof}
By Lemma~\ref{lem:pqpol}, we already know that $p_{n,k}$ and $q_{n,k}$
are polynomials in $A,B,C,D$. The assertion to be shown here is
the \textit{integrality} of the coefficients of these polynomials.
In view of Part~(4) of Remark~\ref{rem:1}, since we also know by
Lemma~\ref{lem:pqpol} that $p_{n,k}$ and $q_{n,k}$ are {\it
  homogeneous} polynomials in $A,B,C,D$, we may without loss of
generality set $B=1$, which we shall do for the rest of the argument.

By inspecting the expressions in \eqref{eq:pnk} and \eqref{eq:qnk},
one sees that the only possible obstacle to integrality of coefficients 
might be powers of $2$ appearing in denominators.

\medskip
The claim is that no denominators occur after one has replaced
$E^2$ by $A^2-4CD$.
The reader should keep in mind that both expressions 
between parentheses in \eqref{eq:pnk} and
\eqref{eq:qnk} are of the form $\Pol(E)-\Pol(-E)$, where $\Pol(E)$ is some
polynomial in $E$ with coefficients depending on $n,k,A,C$. 
(We have again suppressed the dependence
on variables other than $E$ for the sake of better readability. The
reader should recall that $B=1$.) 
Moreover, the polynomial $\Pol(E)$ may be
expanded as a homogeneous polynomial in $E$ and $(A+2C+E)/2$
with integer coefficients. Indeed, that the coefficients are
integers is obvious for $q_{n,k}$. For $p_{n,k}$ this is also
the case since
$$
\binom {k+j}k\binom {n-j}{k-j}\frac { 2 k j}{k+j}=
2j\binom {k+j-1}{k-1}\binom {n-j}{k-j}
$$
and 
$$
\binom {k+j}k\binom {n-j}{k-j}\frac { k- j}{k+j}=
\binom {k+j}k\binom {n-j}{k-j}
-2\binom {k+j-1}k\binom {n-j}{k-j}.
$$
So, what we should examine are differences of the form
\begin{multline} \label{eq:EA}
E^a\left(\frac {A+2C+E} {2}\right)^b
-(-E)^a\left(\frac {A+2C-E} {2}\right)^b\\
=
E^a\left(\left(\frac {A+2C+E} {2}\right)^b
\pm\left(\frac {A+2C-E} {2}\right)^b\right).
\end{multline}
We claim that, after the substitution 
$E^2\to A^2-4CD$, expressions of this form
have integer coefficients.
Due to the sign $\pm$ in the previous expression, we
would have to discuss two different cases. We content ourselves
with treating the case where the sign is positive, the
other case being completely analogous. 
Using the multinomial theorem, we expand
\begin{equation*} 
\left(\frac {A+2C+E} {2}\right)^b
+\left(\frac {A+2C-E} {2}\right)^b
=2^{1-b}\sum_{\ell_1+\ell_2+2\ell_3=b} 
\frac {2^{\ell_2}\,b!} {\ell_1!\,\ell_2!\,(2\ell_3)!}
A^{\ell_1}C^{\ell_2}E^{2\ell_3},
\end{equation*}
and then substitute $A^2-4CD$ for $E^2$. 
Expanding again using the binomial theorem, this leads to the multiple sum
$$
2^{1-b}\sum_{\ell_1+\ell_2+2\ell_3=b} 
\frac {2^{\ell_2}\,b!} {\ell_1!\,\ell_2!\,(2\ell_3)!}
A^{\ell_1}C^{\ell_2}
\sum_{\ell_4=0}^{\ell_3} (-1)^{\ell_3-\ell_4}
\binom {\ell_3}{\ell_4}4^{\ell_3-\ell_4}
A^{2\ell_4}C^{\ell_3-\ell_4}D^{\ell_3-\ell_4}.
$$
We let $a=\ell_1+2\ell_4$, $d=\ell_3-\ell_4$. With this re-indexing,
the above sum is converted into
$$
2^{1-b}\underset{a+d\le b}{\sum_{a,d\ge0} }\sum_{\ell_4\ge0}
\frac {(-1)^d\,2^{b-a}\,b!\,(d+\ell_4)!}
      {(a-2\ell_4)!\,(b-a-2d)!\,(2d+2\ell_4)!\,\ell_4!\,d!}
A^{a}C^{b-a-d}D^d.
$$
The sum over $\ell_4$ can be written in hypergeometric notation,
which leads to the expression
$$
2^{1-b}\underset{a+d\le b}{\sum_{a,d\ge0} }
A^{a}C^{b-a-d}D^d
\frac {(-1)^d\,2^{b-a}\,b!}
      {a!\,(b-a-2d)!\,(2d)!}
{}_2 F_1\!\left[\begin{matrix}
-\frac {a} {2},-\frac {a} {2}+\frac {1} {2}\\d+\frac {1} {2}
\end{matrix}; 1\right].
$$
This $_2F_1$-series can be summed by means of Gau{\ss}' summation
formula 
(see \cite[(1.7.6); Appendix (III.3)]{SlatAC})
\begin{equation*} 
{} _{2} F _{1} \!\left [ \begin{matrix} { a, b}\\ { c}\end{matrix} ;
  {\displaystyle 
   1}\right ]  = \frac {\Gamma ( c)\,\Gamma( c-a-b)} {\Gamma( c-a)\,
  \Gamma( c-b)}.
\end{equation*}
After some simplification, one arrives at
\begin{multline*}
\underset{a+d\le b}{\sum_{a,d\ge0} }
(-1)^dA^{a}C^{b-a-d}D^d\,
\frac {a+2d} {a+d}\binom {a+d}a\binom b{b-a-2d}\\
=\underset{a+d\le b}{\sum_{a,d\ge0} }
(-1)^dA^{a}C^{b-a-d}D^d
\left(\binom {a+d}a
+\binom {a+d-1}a\right)\binom b{b-a-2d}.
\end{multline*}
The integrality of coefficients in expressions of the form
\eqref{eq:EA} is now obvious. As we argued above, this implies
integrality of coefficients in $q_{n,k}$, and, since --- in comparison
to $q_{n,k}$ --- the expression \eqref{eq:pnk} for $p_{n,k}$ contains
an additional $2$ in the denominator, it implies that all
coefficients of $p_{n,k}$ are integers or half-integers (i.e.,
odd integers divided by~2). To see that the latter
coefficients can in fact not be
half-integers, let us suppose for a contradiction that $p_{n,k}$
does contain half-integer coefficients, where $k$ is minimal.
Because of Lemma~\ref{lem:pqpol}, which says 
in particular that $p_{n,0}=q_{n,0}=1$, we must have $k>0$. We then
consider the coefficient of $z^{k}$ on the left-hand side of  
\eqref{eq:idPQ}. It contains the term $p_{n,k}q_{n,0}=p_{n,k}$,
and otherwise only 
terms involving $p_{n,j}$'s with $j<k$ and $q_{n,j}$'s with $j\le k$. Since 
the coefficient of $z^{k}$ on the right-hand side
equals zero, this means that $p_{n,k}$ is given in terms of
expressions which we know by our minimality assumption to be
polynomials in $A,B,C,D$ with integer coefficients. The 
coefficient $p_{n,k}$ must therefore have the same property,
contradicting our assumption that $p_{n,k}$ contains half-integer
coefficients. This completes the proof of the lemma. 
\end{proof}

\section{Free subgroup numbers for lifts of $\PSL_2(\Z)$ modulo
prime powers}
\label{Sec:Freem}

\noindent For integers $m$, $\lambda\geq1$, let $f_\la(m)$ denote the
number of free subgroups  
of index $6m\la$ in the lift $\Gamma_m(3)$ of $\PSL_2(\mathbb Z)$ defined by
\eqref{eq:Gam};
in particular, $f_\la:=f_\la(1)$ is the number of free subgroups in
the inhomogeneous modular group $\PSL_2(\Z)$ 
of index $6\la$. Our original motivation, when embarking on this
project, was to establish the following result, which now is an easy
consequence of Theorem~\ref{thm:ABCD}. 

\begin{theorem} \label{thm:Freem}
Let $m$ be a positive integer,
$p$ a prime number with $p\ge 5,$ and let $\al$ be a positive integer. 
Then the generating function $F_m(z)=1+\sum_{\la\ge1}f_\la(m)z^\la,$
when coefficients are reduced modulo $p^\al,$ can be represented as a
rational function. Equivalently, the sequence $(f_\la(m))_{\la\ge1}$
is ultimately periodic modulo $p^\alpha$.
\end{theorem}

\begin{proof}
It is well-known that $F_m(z)$ satisfies the differential
equation\footnote{Cf.\ \cite[Formula~(8.1)]{KKM} for this, or
  Formula~(18) in \cite[Sec.~2.3]{MuHecke} with $F_m(z) =
  1+z\beta_{\Gamma_m(3)}(z)$ for the same result in a more general
  context.} 
\begin{equation} \label{eq:fmdiff}
(1-(6m-2)z)F_m(z)-6mz^2F_m'(z)-zF_m^2(z)
-1-(1-6m+5m^2)z=0.
\end{equation}
By specialising $A=6m-2$, $B=6m$, $C=1$, and  $D=1-6m+5m^2$ in 
Theorem~\ref{thm:ABCD}, we infer that
there exist polynomials $P^{(m)}_n(z)$ and $Q^{(m)}_n(z)$ such that
$R^{(m)}_n(z)=P^{(m)}_n(z)/Q^{(m)}_n(z)$ satisfies
\begin{multline} \label{eq:Rmdiff}
(1-(6m-2)z)R^{(m)}_n(z)
-6mz^2\left(R^{(m)}_n\right)'(z)-z\left(R^{(m)}_n\right)^2(z)
-1-(1-6m+5m^2)z\\
=
-\frac {5m^{2n+2}z^{2n+1}} {\big(Q^{(m)}_n\big)^2(z)}
\prod _{\ell=1} ^{n}
(6\ell+1)(6\ell+5).
\end{multline}
If we now choose $n$ large enough then, since $p^\al$ is coprime to $6$, 
the right-hand
side will vanish modulo $p^\al$.
Since the differential equation \eqref{eq:fmdiff} determines a 
unique formal power series solution, even when coefficients are
reduced modulo a prime~$p$ (this is seen in the same way as in
the beginning of the proof of Lemma~\ref{lem:pqeind}), 
this proves our claim.
\end{proof}

\section{Some more precise results}
\label{Sec:Periods}

\noindent Recall that the generating function 
$F_m(z)=1+\sum_{\la\ge1}f_\la(m)\,z^\la$ for the numbers of free
subgroups in $\Gamma_m(3)$ satisfies the 
Riccati differential equation \eqref{eq:fmdiff},
which is the special case of \eqref{eq:diffeqABCD} in which $A=6m-2,$
$B=6m,$ $C=1,$ and $D=1-6m+5m^2$.

The following two lemmas are the key for deriving the main result
of this section, Theorem~\ref{thm:Freep}, which identifies
the denominators of the 
rational functions which one obtains when
the coefficients of the generating function
$1+\sum_{\la\ge1}f_\la(m) z^\la$ are reduced modulo a prime power
$p^\al$ with $p\ge5$. 
In the statement of the lemmas, and also later, we write 
$$f(z)=g(z)~\text {modulo}~m$$ 
to mean that the coefficients
of $z^i$ in the power series
$f(z)$ and $g(z)$ agree modulo~$m$ for all $i$.

\begin{lemma} \label{lem:Qnp1}
Let $p$ be a prime with $p\equiv1$~{\em(mod $6$)}.
For $A=6m-2,$ $B=6m,$ $C=1,$ $D=1-6m+5m^2$ {\em(}and, hence, $E=4m${\em),}
and $n\equiv\frac {p-1} {6},\frac {5(p-1)} {6}$~{\em(mod~$p$),} we have
$$
Q_n(z)=Q_{(p-1)/6}(z)\quad \text{\em modulo $p$},
$$
with $Q_n(z)$ as given in the assertion of Theorem~{\em\ref{thm:ABCD}}.
Furthermore, if $p\nmid m$, then 
the polynomial $Q_{(p-1)/6}(z)$ has degree $(p-1)/6$ in~$z$
when coefficients are reduced modulo~$p$, whereas 
$Q_{(p-1)/6}(z)=1$~{\em modulo}~$p$ if $p\mid m$.
\end{lemma}

\begin{proof}
We specialise the result of Theorem~\ref{thm:ABCD} to $A=6m-2$,
$B=6m$, $C=1$, $D=1-6m+5m^2$, to see that $R_n(z)=P_n(z)/Q_n(z)$ solves
\eqref{eq:fmdiff}, where $Q_n(z)=\sum_{k=0}^n q_{n,k}z^k$,
with coefficients given by
\begin{multline} \label{eq:qnkFree}
q_{n,n-k}=\frac {(-1)^n6^{n-k}m^{n-k}} {(\frac {2} {3}-k)_{2k+1}}
\Bigg( \left(\tfrac {5} {6}\right)_{n+1}
\sum_{j=0}^k \binom {k+j}k\binom {n-j}{k-j}
\left(-\tfrac 23
  + j + 1\right)_{ 
    k - j} \left(\tfrac {1} {6}\right)_j\\
- \left(\tfrac {1} {6}\right)_{n+1}
\sum_{j=0}^k \binom {k+j}k\binom {n-j}{k-j}
\left(\tfrac 23
  + j + 1\right)_{ 
    k - j} \left(\tfrac {5} {6}\right)_j
\Bigg).
\end{multline}
The reader should observe that the only dependence on $m$ of
$q_{n,n-k}$ is in the term $m^{n-k}$ which occurs right in the front
of the above expression. Hence, since
by Theorem~\ref{thm:ABCD} we know that $q_{n,0}=1$, we have
$Q_{(p-1)/6}(z)=1$~modulo~$p$ if $p\mid m$. On the other hand, 
if $p\nmid m$, from \eqref{eq:qnn} it is easy to see that
$q_{(p-1)/6,(p-1)/6}\hbox{${}\not\equiv{}$}0$~(mod~$p$)
for $A=6m-2$, $B=6m$, $C=1$, $D=1-6m+5m^2$, $E=4m$.
(As a matter of fact, with this choice of parameters, we have
$\Pi_+\hbox{${}\not\equiv{}$}0 $~(mod~$p$), whereas
$\Pi_-\equiv0 $~(mod~$p$).)
So, indeed, if $p\nmid m$ 
the polynomial $Q_{(p-1)/6}(z)$ has degree $(p-1)/6$
in~$z$ when coefficients are reduced modulo~$p$.

What remains to show is that all coefficients $q_{n,n-k}$,
$0\le k\le n-\frac {p-1} {6}-1$, 
are divisible by $p$ for the specialisation of the parameters
$A,B,C,D,E$ considered here.

In order to see this, we start by
writing the sums over $j$ in \eqref{eq:qnkFree} 
in hypergeometric notation, thus obtaining
\begin{multline*} 
q_{n,n-k}=\frac {(-1)^n6^{n-k}m^{n-k}} {(\frac {2} {3}-k)_{2k+1}}
\Bigg( \left(\tfrac {5} {6}\right)_{n+1}
\binom {n}{k}
\left(\tfrac 13
   \right)_{ 
    k }
\,{}_3 F_2\!\left[\begin{matrix}
\frac {1} {6},k+1,-k\\\frac {1} {3},-n
\end{matrix};1\right]
\\
- \left(\tfrac {1} {6}\right)_{n+1}
\binom {n}{k}
\left(\tfrac 53
  \right)_{ 
    k } 
\,{}_3 F_2\!\left[\begin{matrix}
\frac {5} {6},k+1,-k\\\frac {5} {3},-n
\end{matrix};1\right]
\Bigg).
\end{multline*}
Next we apply the transformation formula \eqref{eq:T3204} one more
time, here with $a=-k$.
This converts the last expression into
\begin{align}
\notag
q_{n,n-k}&=\frac {(-1)^n6^{n-k}m^{n-k}} {(\frac {2} {3}-k)_{2k+1}}
\Bigg( \left(\tfrac {5} {6}\right)_{n+1}
\frac {(-1)^k} {k!}
\left(\tfrac 13
   \right)_{ 
    k }\,
\left(-n-k-\tfrac {5} {6}\right)_k
\,{}_3 F_2\!\left[\begin{matrix}\vphantom{\Big(}
\frac {1} {6},-k-\frac {2} {3},-k\\\frac {1} {3},
-n-k-\frac {5} {6}
\end{matrix};1\right]
\\
\notag
&\kern3cm
- \left(\tfrac {1} {6}\right)_{n+1}
\frac {(-1)^k} {k!}
\left(\tfrac 53
  \right)_{ 
    k } 
\left(-n-k-\tfrac {1} {6}\right)_k
\,{}_3 F_2\!\left[\begin{matrix}\vphantom{\Big(}
\frac {5} {6},-k+\frac {2} {3},-k\\\frac {5} {3},
-n-k-\frac {1} {6}
\end{matrix};1\right]
\Bigg)\\
\notag
&
={(-1)^n6^{n-k}m^{n-k}}
\Bigg(
\sum_{j=0}^k
\frac {(-1)^{k+j}} {k!}\binom kj
\frac {\left(\tfrac {5} {6}-j\right)_{n+k+1}} 
{\left(\tfrac {2} {3}-j\right)_{k+1}} \\
&\kern7cm
+
\sum_{j=0}^k
\frac {(-1)^{k+j}} {k!}\binom kj
\frac {\left(\tfrac {1} {6}-j\right)_{n+k+1}} 
{\left(-\tfrac {2} {3}-j\right)_{k+1}} 
\Bigg).
\label{eq:qnkp}
\end{align}
We claim that, if $n\equiv\frac {p-1} {6},\frac {5(p-1)} {6}$~(mod~$p$) 
and $n\ge k+\frac {p-1} {6}+1$, then
each summand in the sums over $j$ is divisible by $p$.

To see this, we first consider the expression
$$
\frac {(-1)^{k+j}} {k!}\binom kj
\frac {\left(\tfrac {5} {6}-j\right)_{n+k+1}} 
{\left(\tfrac {2} {3}-j\right)_{k+1}} .
$$
Let $v_p(\al)$ denote the $p$-adic valuation of the
integer $\al$, that is, the maximal exponent $e$ such that $p^e$
divides $\al$. By variations of the well-known formula of
Legendre \cite[p.~10]{LegeAA} for the $p$-adic valuation of
factorials, we have
\begin{multline} \label{eq:expp}
v_p\Bigg(
\frac {(-1)^{k+j}} {k!}\binom kj
\frac {\left(\tfrac {5} {6}-j\right)_{n+k+1}} 
{\left(\tfrac {2} {3}-j\right)_{k+1}} 
\Bigg)
=
\sum_{\ell=1}^{\infty}
\Bigg(
-\fl{\frac {j} {p^\ell}}
-\fl{\frac {k-j} {p^\ell}}
+\fl{\frac {n+k-j+\frac {1} {6}(p^\ell+5)} {p^\ell}}\\
-\fl{\frac {-j+\frac {1} {6}(p^\ell-1)} {p^\ell}}
-\fl{\frac {k-j+\frac {1} {3}(p^\ell+2)} {p^\ell}}
+\fl{\frac {-j+\frac {1} {3}(p^\ell-1)} {p^\ell}}
\Bigg).
\end{multline}
We need to prove that this sum over $\ell$ is at least $1$.
In the case where $\ell=1$, the summand reduces to
\begin{multline*}
-\fl{\frac {j} {p}}
-\fl{\frac {k-j} {p}}\\
+\fl{\frac {n+k-j+\frac {5} {6}} {p}+\frac {1} {6}}
-\fl{\frac {-j-\frac {1} {6}} {p}+\frac {1} {6}}
-\fl{\frac {k-j+\frac {2} {3}} {p}+\frac {1} {3}}
+\fl{\frac {-j-\frac {1} {3}} {p}+\frac {1} {3}}.
\end{multline*}
Let us write
$N=\{n/p\}$,
$K=\{k/p\}$, and
$J=\{j/p\}$,
where $\{\al\}:=\al-\fl{\al}$ denotes the fractional
part of $\al$. With this notation, the last expression becomes
\begin{multline} \label{eq:vp} 
\fl{\frac {n} {p}}
-\fl{\frac {k} {p}}
-\fl{K-J}\\
+\fl{N+K-J+\frac {5} {6p}+\frac {1} {6}}
-\fl{-J-\frac {1} {6p}+\frac {1} {6}}
-\fl{K-J+\frac {2} {3p}+\frac {1} {3}}
+\fl{-J-\frac {1} {3p}+\frac {1} {3}}.
\end{multline}
If $n\equiv\frac {p-1} {6}$~(mod~$p$), then $N=\frac {p-1}{6p}$. Hence,
$$\fl{N+K-J+\frac {5} {6p}+\frac {1} {6}}
-\fl{K-J+\frac {2} {3p}+\frac {1} {3}}=0.
$$
Moreover, since $k\le n-\frac {p-1} {6}-1$, we have $\fl{\frac {n} {p}}
-\fl{\frac {k} {p}}
\ge1$. Since, trivially, $-\fl{K-J}\ge0$ and
$
-\fl{-J-\frac {1} {6p}+\frac {1} {6}}
+\fl{-J-\frac {1} {3p}+\frac {1} {3}}\ge0
$,
the expression in \eqref{eq:vp} is at least~$1$.

On the other hand, if $n\equiv\frac {5(p-1)} {6}$~(mod~$p$), then 
$N=\frac {5(p-1)} {6p}$. 
In this case, we have 
$$-\fl{K-J}
+\fl{N+K-J+\frac {5} {6p}+\frac {1} {6}}
=1,$$ 
we have certainly $\fl{\frac {n} {p}}
-\fl{\frac {k} {p}}
\ge0$,
and, again, $
-\fl{-J-\frac {1} {6p}+\frac {1} {6}}
+\fl{-J-\frac {1} {3p}+\frac {1} {3}}\ge0
$. If we suppose that $-\fl{K-J+\frac {2} {3p}+\frac {1} {3}}=-1$, then
$K$ must be at least $\frac {2} {3p}(p-1)$, or, equivalently, 
$k\equiv\frac {2} {3}(p-1),\frac {2} {3}(p-1)+1,\dots,p-1$~(mod~$p$). 
But, because of $n\ge k+\frac {p-1} {6}+1$ and 
$n\equiv\frac {5(p-1)} {6}$~(mod~$p$), this implies 
$\fl{\frac {n} {p}}
-\fl{\frac {k} {p}}
\ge1$. Hence, in all cases the expression in \eqref{eq:vp} is at least~1.

Now we have to discuss the summand in \eqref{eq:expp} when $\ell\ge2$.
Similarly to above, we write
$N=\{n/p^\ell\}$,
$K=\{k/p^\ell\}$, and
$J=\{j/p^\ell\}$. Using this notation, we may rewrite this summand as
\begin{multline} \label{eq:vpell}
\fl{\frac {n} {p^\ell}}
-\fl{\frac {k} {p^\ell}}
-\fl{K-J}
+\fl{N+K-J+\frac {1} {6}+
\frac {5} {6\cdot p^\ell}}
-\fl{-J+\frac {1} {6}-\frac {1} {6\cdot p^\ell}}\\
-\fl{K-J+\frac {1} {3}+\frac {2} {3\cdot p^\ell}}
+\fl{-J+\frac {1} {3}-\frac {1} {3\cdot p^\ell}}.
\end{multline}
We have $-\fl{-J+\frac {1} {6}-\frac {1} {6\cdot p^\ell}}
+\fl{-J+\frac {1} {3}-\frac {1} {3\cdot p^\ell}}\ge0$,
$-\fl{K-J}\ge0$, and, since
$n\ge k+\frac {p-1} {6}+1$, also $\fl{\frac {n} {p^\ell}}
-\fl{\frac {k} {p^\ell}}\ge0$.
If 
\begin{equation} \label{eq:n+k-j} 
\fl{N+K-J+\frac {1} {6}+
\frac {5} {6\cdot p^\ell}}
-\fl{K-J+\frac {1} {3}+\frac {2} {3\cdot p^\ell}}\ge0,
\end{equation}
then \eqref{eq:vpell} is non-negative. There are two cases 
where \eqref{eq:n+k-j} is violated. It should be noted that then
the left-hand side of \eqref{eq:n+k-j} equals $-1$.
Violation of \eqref{eq:n+k-j} can either happen if 
$$
N+K-J+\frac {1} {6}+
\frac {5} {6\cdot p^\ell}
<0\le K-J+\frac {1} {3}+\frac {2} {3\cdot p^\ell},
$$
in which case necessarily $K-J<0$ and thus actually $-\fl{K-J}=1$, or if
$$
N+K-J+\frac {1} {6}+
\frac {5} {6\cdot p^\ell}
<1\le K-J+\frac {1} {3}+\frac {2} {3\cdot p^\ell}.
$$
In the latter case, 
we have in particular $K\ge \frac {2} {3}-\frac {2}
  {3\cdot p^\ell}$, and on the other hand
$$
N
\le
N+K-J-\frac {2} {3}+
\frac {2} {3\cdot p^\ell}<
\frac {1} {6}
-\frac {1} {6\cdot p^\ell}.
$$
Again, since $n\ge k+\frac {p-1} {6}+1$, 
this implies that actually $\fl{\frac {n} {p^\ell}}
-\fl{\frac {k} {p^\ell}}\ge1$, so that in all cases the expression in
\eqref{eq:vpell} is non-negative.

This proves that the sum on the right-hand side of \eqref{eq:expp} is
at least $1$, which, in turn, means that all summands in the
first sum over~$j$ on the right-hand side of \eqref{eq:qnkp} are
divisible by~$p$.

Next, we consider the expression
$$
\frac {(-1)^{k+j}} {k!}\binom kj
\frac {\left(\tfrac {1} {6}-j\right)_{n+k+1}} 
{\left(-\tfrac {2} {3}-j\right)_{k+1}} .
$$
Here, we see that
\begin{multline} \label{eq:expp2}
v_p\Bigg(
\frac {(-1)^{k+j}} {k!}\binom kj
\frac {\left(\tfrac {1} {6}-j\right)_{n+k+1}} 
{\left(-\tfrac {2} {3}-j\right)_{k+1}} 
\Bigg)\\
=
\sum_{\ell=1}^{\infty}
\Bigg(
-\fl{\frac {j} {p^\ell}}
-\fl{\frac {k-j} {p^\ell}}
+\fl{\frac {n+k-j+\frac {1} {6}(5\cdot p^\ell+1)} {p^\ell}}
\kern3cm\\
-\fl{\frac {-j+\frac {5} {6}(p^\ell-1)} {p^\ell}}
-\fl{\frac {k-j+\frac {2} {3}(p^\ell-1)} {p^\ell}}
+\fl{\frac {-j+\frac {1} {3}(2\cdot p^\ell-5)} {p^\ell}}
\Bigg).
\end{multline}
Again, we have to prove that this sum over $\ell$ is at least $1$.
In the case where $\ell=1$, the summand reduces to
\begin{multline*}
-\fl{\frac {j} {p}}
-\fl{\frac {k-j} {p}}\\
+\fl{\frac {n+k-j+\frac {1} {6}} {p}+\frac {5} {6}}
-\fl{\frac {-j-\frac {5} {6}} {p}+\frac {5} {6}}
-\fl{\frac {k-j-\frac {2} {3}} {p}+\frac {2} {3}}
+\fl{\frac {-j-\frac {5} {3}} {p}+\frac {2} {3}}.
\end{multline*}
With the notation
$N=\{n/p\}$,
$K=\{k/p\}$, and
$J=\{j/p\}$, as before,
the last expression becomes
\begin{multline} \label{eq:vp2} 
\fl{\frac {n} {p}}
-\fl{\frac {k} {p}}
-\fl{K-J}\\
+\fl{N+K-J+\frac {1} {6p}+\frac {5} {6}}
-\fl{-J-\frac {5} {6p}+\frac {5} {6}}
-\fl{K-J-\frac {2} {3p}+\frac {2} {3}}
+\fl{-J-\frac {5} {3p}+\frac {2} {3}}.
\end{multline}
If $n\equiv\frac {p-1} {6}$~(mod~$p$), then $N=\frac {p-1} {6p}$. Hence,
$$\fl{N+K-J+\frac {1} {6p}+\frac {5} {6}}
-\fl{K-J}=1.
$$
Moreover, since $k\le n-\frac {p-1} {6}-1$, we have $\fl{\frac {n} {p}}
-\fl{\frac {k} {p}}
\ge1$. If
$$
-\fl{-J-\frac {5} {6p}+\frac {5} {6}}
+\fl{-J-\frac {5} {3p}+\frac {2} {3}}=-1,
$$
then $J=\frac {2(p-1)} {3p},\frac {2(p-1)} {3p}+1,\dots,p-1$, 
and in all cases it follows that
$$
-\fl{K-J-\frac {2} {3p}+\frac {2} {3}}\ge0.
$$
Therefore the expression in \eqref{eq:vp2} is at least~$1$.

On the other hand, if $n\equiv\frac {5(p-1)} {6}$~(mod~$p$), then 
$N=\frac {5(p-1)} {6p}$. 
In this case, we have 
$$
\fl{N+K-J+\frac {1} {6p}+\frac {5} {6}}
-\fl{K-J-\frac {2} {3p}+\frac {2} {3}}=1$$
and $\fl{\frac {n} {p}}
-\fl{\frac {k} {p}} \ge0$. Again, if 
$
-\fl{-J-\frac {5} {6p}+\frac {5} {6}}
+\fl{-J-\frac {5} {3p}+\frac {2} {3}}=-1,
$
then $J=\frac {2(p-1)} {3p},\frac {2(p-1)} {3p}+1,\dots,p-1$, 
and in all cases it follows that
$$-\fl{K-J}\ge -\fl{K-\frac {2(p-1)} {3p}}\ge0.$$
If $\fl{K-J}=0$ (that is, if $-\fl{K-J}\ne1$), 
so that in particular $K\ge\frac {2(p-1)} {3p}$, then,
since $n\equiv\frac {5(p-1)} {6}$~(mod~$p$) and $n\ge k+\frac {p-1}
{6}+1$, we actually have $\fl{\frac {n} {p}}-\fl{\frac {k} {p}} \ge1$.
Hence, in all cases the expression in \eqref{eq:vp2} is at least~1.

\medskip
We turn to the summand in \eqref{eq:expp2} when $\ell\ge2$.
Writing
$N=\{n/p^\ell\}$,
$K=\{k/p^\ell\}$, and
$J=\{j/p^\ell\}$, this summand acquires the form
\begin{multline} \label{eq:vpell2}
\fl{\frac {n} {p^\ell}}
-\fl{\frac {k} {p^\ell}}
-\fl{K-J}
+\fl{N+K-J+\frac {5} {6}+\frac {1}
  {6\cdot p^\ell}}\\
-\fl{-J+\frac {5} {6}-\frac {5} {6\cdot p^\ell}}
-\fl{K-J+\frac {2} {3}-\frac {2} {3\cdot p^\ell}}
+\fl{-J+\frac {2} {3}-\frac {5}
  {3\cdot p^\ell}}.
\end{multline}
Clearly, 
we have
$-\fl{K-J}\ge0$, we have
$$\fl{N+K-J+\frac {5} {6}+\frac {1}
  {6\cdot p^\ell}}
-\fl{K-J+\frac {2} {3}-\frac {2} {3\cdot p^\ell}}\ge0,$$ 
and, since
$n\ge k+\frac {p-1} {6}+1$, also $\fl{\frac {n} {p^\ell}}
-\fl{\frac {k} {p^\ell}}\ge0$.
If 
\begin{equation} \label{eq:-j} 
-\fl{-J+\frac {5} {6}-\frac {5} {6\cdot p^\ell}}
+\fl{-J+\frac {2} {3}-\frac {5}
  {3\cdot p^\ell}}\ge0,
\end{equation}
then \eqref{eq:vpell2} is non-negative. The inequality \eqref{eq:-j}
is only violated if $\frac {5} {6}-\frac {5} {6\cdot p^\ell}\ge
J>\frac {2} {3}-\frac {5}
  {3\cdot p^\ell}$, in which case the
left-hand side of \eqref{eq:-j} equals $-1$.
It follows that $K-J+\frac {2} {3}-\frac {2} {3\cdot p^\ell}<K+\frac
{1} {p^\ell}\le 1$, so that $\fl{K-J+\frac {2} {3}-\frac {2} {3\cdot
    p^\ell}}\le0$, that is, equivalently, 
$-\fl{K-J+\frac {2} {3}-\frac {2} {3\cdot
    p^\ell}}\ge0$.
If we now suppose that both $-\fl{K-J}$ and
$\big\lfloor N+K-J+\frac {5} {6}+\frac {1}
  {6\cdot p^\ell}\big\rfloor$ are zero (that is, not~$1$), then $K\ge J$
and $N+K-J+\frac {5} {6}+\frac {1}
  {6\cdot p^\ell}<1$, and this entails the inequality chain
$$
1>N+K-J+\frac {5} {6}+\frac {1}
  {6\cdot p^\ell}\ge N+K+\frac {1} {p^\ell}.
$$
On the other hand,
since $K\ge J> \frac {2} {3}-\frac {5} {3\cdot p^\ell}$, the number
$N$ can be at most $1/3$. 
Because of $n\ge k+\frac {p-1} {6}+1$, 
this implies that actually $\fl{\frac {n} {p^\ell}}
-\fl{\frac {k} {p^\ell}}\ge1$, so that in all cases the expression in
\eqref{eq:vpell2} is non-negative.

This proves that the sum on the right-hand side of \eqref{eq:expp2} is
at least $1$, which, in its turn, means that all summands in the
second sum over~$j$ on the right-hand side of \eqref{eq:qnkp} are
divisible by~$p$, completing the proof of the lemma. 
\end{proof}

\begin{lemma} \label{lem:Qnp2}
Let $p$ be a prime with $p\equiv5$~{\em(mod $6$)}.
For $A=6m-2,$ $B=6m,$ $C=1,$ $D=1-6m+5m^2$ {\em(}and, hence, $E=4m${\em),}
and $n\equiv\frac {p-5} {6},\frac {5p-1} {6}$~{\em(mod~$p$),} we have
$$
Q_n(z)=Q_{(p-5)/6}(z)\quad \text{\em modulo $p$},
$$
with $Q_n(z)$ as given in Theorem~{\em\ref{thm:ABCD}}.
Furthermore, if $p\nmid m$, then 
the polynomial\break $Q_{(p-5)/6}(z)$ has degree $(p-5)/6$
in~$z$ when coefficients are reduced modulo~$p$, whereas 
$Q_{(p-5)/6}(z)=1$~{\em modulo}~$p$ if $p\mid m$.
\end{lemma}

\begin{proof}
This assertion may be established in a fashion similar to the proof of
Lemma~\ref{lem:Qnp1}. 
The assertion about the degree of $Q_{(p-5)/6}(z)$
 is straightforward to see.

In order to prove the congruence, we would again use the form
\eqref{eq:qnkp} for $q_{n,n-k}$.
Here, we have
\begin{multline} \label{eq:expp3}
v_p\Bigg(
\frac {(-1)^{k+j}} {k!}\binom kj
\frac {\left(\tfrac {5} {6}-j\right)_{n+k+1}} 
{\left(\tfrac {2} {3}-j\right)_{k+1}} 
\Bigg)
=
\underset{\ell\text{ even}}{\sum_{\ell=1}^{\infty}}
\Bigg(
-\fl{\frac {j} {p^\ell}}
-\fl{\frac {k-j} {p^\ell}}
+\fl{\frac {n+k-j+\frac {1} {6}(p^\ell+5)} {p^\ell}}\\
-\fl{\frac {-j+\frac {1} {6}(p^\ell-1)} {p^\ell}}
-\fl{\frac {k-j+\frac {1} {3}(p^\ell+2)} {p^\ell}}
+\fl{\frac {-j+\frac {1} {3}(p^\ell-1)} {p^\ell}}
\Bigg)\\
+
\underset{\ell\text{ odd}}{\sum_{\ell=1}^{\infty}}
\Bigg(
-\fl{\frac {j} {p^\ell}}
-\fl{\frac {k-j} {p^\ell}}
+\fl{\frac {n+k-j+\frac {5} {6}(p^\ell+1)} {p^\ell}}
\kern4cm\\
-\fl{\frac {-j+\frac {1} {6}(5p^\ell-1)} {p^\ell}}
-\fl{\frac {k-j+\frac {2} {3}(p^\ell+1)} {p^\ell}}
+\fl{\frac {-j+\frac {1} {3}(2p^\ell-1)} {p^\ell}}
\Bigg)
\end{multline}
and
\begin{multline} \label{eq:expp4}
v_p\Bigg(
\frac {(-1)^{k+j}} {k!}\binom kj
\frac {\left(\tfrac {1} {6}-j\right)_{n+k+1}} 
{\left(-\tfrac {2} {3}-j\right)_{k+1}} 
\Bigg)\\=
\underset{\ell\text{ even}}{\sum_{\ell=1}^{\infty}}
\Bigg(
-\fl{\frac {j} {p^\ell}}
-\fl{\frac {k-j} {p^\ell}}
+\fl{\frac {n+k-j+\frac {1} {6}(5p^\ell+1)} {p^\ell}}\\
-\fl{\frac {-j+\frac {5} {6}(p^\ell-1)} {p^\ell}}
-\fl{\frac {k-j+\frac {2} {3}(p^\ell-1)} {p^\ell}}
+\fl{\frac {-j+\frac {1} {3}(2p^\ell-5)} {p^\ell}}
\Bigg)\\
+
\underset{\ell\text{ odd}}{\sum_{\ell=1}^{\infty}}
\Bigg(
-\fl{\frac {j} {p^\ell}}
-\fl{\frac {k-j} {p^\ell}}
+\fl{\frac {n+k-j+\frac {1} {6}(p^\ell+1)} {p^\ell}}
\kern4cm\\
-\fl{\frac {-j+\frac {1} {6}(p^\ell-5)} {p^\ell}}
-\fl{\frac {k-j+\frac {1} {3}(2p^\ell-1)} {p^\ell}}
+\fl{\frac {-j+\frac {1} {3}(2p^\ell-4)} {p^\ell}}
\Bigg),
\end{multline}
and these are the expressions which have to be analysed $p$-adically.
We leave it to the interested reader to fill in the details.
\end{proof}

\begin{theorem} \label{thm:Freep}
Let $\al$ be a positive integer.
If $p$ is a prime with $p\equiv1$~{\em(mod $6$),} 
then the generating function $F_m(z)=1+\sum_{\la\ge1}f_\la(m)\,z^\la,$ when 
coefficients are reduced modulo~$p^\al,$ equals 
$\overline P_\al(z)/Q_{(p-1)/6}^\al(z),$
where $\overline P_\al(z)$ is a polynomial in $z$ over the integers,
and $Q_{(p-1)/6}(z)$ is the polynomial of degree $(p-1)/6$
given by the formula in Theorem~{\em\ref{thm:ABCD}} with
$A=6m-2,$ $B=6m,$ $C=1,$ $D=1-6m+5m^2$ and $E=4m$.

On the other hand,
if $p$ is a prime with $p\equiv5$~{\em(mod $6$),} 
then the generating\break function $F_m(z)=1+\sum_{\la\ge1}f_\la(m)\,z^\la,$ when 
coefficients are reduced modulo~$p^\al,$ equals\break
$\widehat P_\al(z)/Q_{(p-5)/6}^\al(z),$
where $\widehat P_\al(z)$ is a polynomial in $z$ over the integers,
and $Q_{(p-5)/6}(z)$ is the polynomial of degree $(p-5)/6$
given by the formula in Theorem~{\em\ref{thm:ABCD}} with
$A=6m-2,$ $B=6m,$ $C=1,$ $D=1-6m+5m^2$ and $E=4m$.
\end{theorem}

\begin{proof}
We proceed by induction on $\al$.

For the start of the induction, we just have to use
Theorem~\ref{thm:ABCD} with 
$A=6m-2$, $B=6m$, $C=1$, $D=1-6m+5m^2$, $E=4m$, and choose
$n=(p-1)/6$ and $n=(p-5)/6$, respectively, as this makes
the product on the right-hand side of \eqref{eq:idABCD} vanish.
By uniqueness of
the solution of \eqref{eq:fmdiff}, it follows that
$$
F_m(z)=\frac {P_{d}(z)} {Q_d(z)}\quad \text {modulo }p,
$$
where $d=(p-1)/6$ or $d=(p-5)/6$, depending on the congruence class of
$p$ modulo~$6$.

Let us now suppose that we have already found an integer polynomial 
$A_\al(z)$ such that
$F_\al(z):=A_\al(z)/Q^\al_d(z)$ 
solves \eqref{eq:fmdiff} modulo~$p^\al$ for some $\al\ge1$.
We choose $n\equiv\frac {p-1} {6}$~(mod~$p$) or
$n\equiv\frac {p-5} {6}$~(mod~$p$), respectively, 
large enough such that the product over
$\ell$ on the right-hand side of \eqref{eq:Rmdiff} 
vanishes modulo~$p^{\al+1}$. 
(The reader should recall that the specialisation 
$A=6m-2$, $B=6m$, $C=1$, $D=1-6m+5m^2$, $E=4m$,
which is considered here, corresponds to the
specialisation of Theorem~\ref{thm:ABCD} considered in
Theorem~\ref{thm:Freem}).
By uniqueness of solution of the
differential equation \eqref{eq:Rmdiff} (see Theorem~\ref{thm:ABCD}),
we must then have 
\begin{equation} \label{eq:RnF}
R_n(z)=F_\al(z)\quad 
\text{modulo }p^\al,
\end{equation}
where $R_n(z)=R^{(1)}_n(z)=P_n(z)/Q_n(z)$ is the solution to
\eqref{eq:Rmdiff} given by Theorem~\ref{thm:ABCD}. 
Consequently, if we consider the difference,
$$
R_n(z)-F_\al(z)=\frac {P_n(z)} {Q_n(z)}-\frac {A_\al(z)} {Q^\al_d(z)}
=
\frac {P_n(z)Q^\al_d(z)-A_\al(z)Q_n(z)} {Q^\al_d(z) Q_n(z)},
$$ 
then, since $Q^\al_d(z) Q_n(z)$ has constant coefficient $1$
(see Theorem~\ref{thm:ABCD}), we know by \eqref{eq:RnF} 
that all coefficients of
the integer polynomial in the numerator of the last fraction must be
divisible by $p^\al$. In other words, there is an integer polynomial
$B_\al(z)$ such that
$$
R_n(z)=F_\al(z)+p^\al \frac {B_\al(z)} {Q^\al_d(z) Q_n(z)}.
$$
If we consider both sides of this equation modulo $p^{\al+1}$,
then the fraction which is multiplied by $p^\al$ may be reduced
modulo~$p$. By Lemmas~\ref{lem:Qnp1} and \ref{lem:Qnp2}, 
this leads to the congruence
\begin{equation} \label{eq:induction} 
R_n(z)=\frac {A_\al(z)} {Q^\al_d(z)}
+p^\al \frac {B_\al(z)} {Q^{\al+1}_d(z)}
\quad \text{modulo }p^{\al+1}.
\end{equation}
So, indeed, the solution $R_n(z)$ to the Riccati differential equation
\eqref{eq:fmdiff}, when coefficients of series are reduced
modulo~$p^{\al+1}$, can be expressed as a rational function with
denominator $Q^{\al+1}_d(z)$, where $d=(p-1)/6$ in case
$p\equiv1$~(mod~$6$), and $d=(p-5)/6$ if
$p\equiv5$~(mod~$6$). This establishes the theorem.
\end{proof}

In what follows, we concentrate on the case where $m=1$, that is, on
the number $f_\la=f_\la(1)$ of free subgroups of index $6\la$ in 
the modular group
$\PSL_2(\Z)$. Note that this means that we consider the case of
Theorems~\ref{thm:ABCD} and \ref{thm:Freep} where
$A=E=4$, $B=6$, $C=1$, and $D=0$.
Our goal is to make Theorem~\ref{thm:Freep} in this special case 
more explicit for the first few prime numbers,  
and to determine precisely the (minimal) period length of the
sequence $(f_\la)_{\la\ge1}$ when it is reduced modulo these prime
numbers. 

\begin{theorem} \label{thm:Free7}
Let $\al$ be a positive integer.
Then the generating function $F(z)=1+\sum_{\la\ge1}f_\la\,z^\la$, when 
coefficients are reduced modulo~$7^\al$, equals $P_\al(z)/(1+2z)^\al$,
where $P_\al(z)$ is a polynomial in $z$ over the integers. In particular,
the sequence $(f_\la)_{\la\ge1}$ is ultimately periodic, with minimal
period equal to $6\cdot 7^{\al-1}$.  
\end{theorem}

\begin{proof}
In view of Theorem~\ref{thm:Freep}, and since $Q_1(z)=1+2z$~modulo~$7$
for $A=E=4$, $B=6$, $C=1$, the only assertion which still
requires proof is the one about the period length.  
A careful examination of the previous arguments
(see in particular \eqref{eq:induction}) reveals that
we have actually shown the stronger statement that
the solution $F(z)$ to \eqref{eq:fmdiff} with $m=1$, 
when coefficients are
reduced modulo~$7^\al$, can be written in the form
$$
F(z)=\text{Pol}_\al(z)+\sum_{k=1}^{\al}f_{-k}7^{k-1}(1+2z)^{-k}
\quad \text{modulo }7^\al,
$$
where $\text{Pol}_\al(z)$ is a polynomial in $z$ with integer
coefficients, and the $f_{-k}$'s are integers. 
Now, we have
$$
7^{k-1}(1+2z)^{-k}=\sum_{\ell\ge0}(-2)^\ell7^{k-1}\binom {k+\ell-1}{k-1}z^\ell.
$$
The term $(-2)^\ell$, $\ell=0,1,\dots$, is periodic modulo~$7^\al$
with (minimal) period 
length $\ph(7^\al)=6\cdot 7^{\al-1}$, where $\ph$ denotes Euler's 
totient function. Indeed, this holds
for $\alpha=1,2$ by inspection, and thus for every $\alpha\geq1$ by a
theorem of Gau{\ss} concerning existence of primitive roots
(cf.\ e.g.\ \cite[Ch.~2, Sec.~5]{BundAA} or 
\cite[pp.~285--287]{Ore}). On the other hand, the period length of 
$$
7^{k-1}\binom {k+\ell-1}{k-1}=
7^{k-1-v_7((k-1)!)}\frac {(k+\ell-1)\cdots(\ell+2)(\ell+1)} 
{(k-1)!7^{-v_7((k-1)!)}},
\quad \ell=0,1,\dots,
$$
when considered modulo $7^\alpha$, 
divides $7^{\al-1}$ for $k\ge2$, since 
$$
k-1-v_7((k-1)!)\ge1
$$
for $k\ge2$. In total, 
$$
(-2)^\ell7^{k-1}\binom {k+\ell-1}{k-1},
\quad \ell=0,1,\dots,
$$
when considered modulo $7^\al$, 
is periodic with period length (exactly) $6\cdot 7^{\al-1}$,
which implies the corresponding assertion of the theorem.
\end{proof}

We have implemented the algorithm which is implicit in the proof of
Theorem~\ref{thm:Freep} (see the Note at the end of the
Introduction). As an illustration, the next
theorem contains the result for the modulus $7^5=16807$.

{
\allowdisplaybreaks
\begin{theorem} \label{thm:Free7^5}
We have
\begin{multline} \label{eq:free7^5}
1+\sum_{\la\ge1}f_\la\,z^\la\\=
4802 z^{25}+9604 z^{23}+14406
   z^{22}+2401 z^{21}+2401
   z^{20}+4802 z^{19}+9947
   z^{18}+9604 z^{17}\\
+10290
   z^{16}+9947 z^{15}+10976
   z^{14}+16464 z^{13}+12691
   z^{12}+2940 z^{11}+8918
   z^{10}\\
+15484 z^9+8722
   z^8+4214 z^7+10829 z^6+6174
   z^5+406 z^4+14896 z^3+11102
   z^2\\
+14168 z+7
+\frac{16451}{1+2z}
+\frac{9562}{(1+2z)^2}
+\frac{2450}{(1+2z)^3}
+\frac{2744}{(1+2z)^4}
+\frac{2401}{(1+2z)^5}\\
\quad \text{\em modulo }7^5.
\end{multline}
\end{theorem}}

The case of powers of $11$ can be handled in an analogous way.
We content ourselves with stating the corresponding result, leaving
the proof to the interested reader. 

\begin{theorem} \label{thm:Free11}
Let $\al$ be a positive integer.
Then the generating function $F(z)=1+\sum_{\la\ge1}f_\la\,z^\la$, when 
coefficients are reduced modulo~$11^\al$, equals $P_\al(z)/(1-z)^\al$,
where $P_\al(z)$ is a polynomial in $z$ over the integers. Moreover,
the sequence $(f_\la)_{\la\ge1}$ is ultimately periodic, with minimal
period length equal to $11^{\al-1}$.  
\end{theorem}

Again, we have implemented the algorithm which leads 
to the above theorem (see the Note at the end of the
Introduction). As an illustration, the next
theorem contains the result for the modulus $11^5=161051$.

{
\allowdisplaybreaks
\begin{theorem} \label{thm:Free11^5}
We have
\begin{multline} \label{eq:free11^5}
1+\sum_{\la\ge1}f_\la\,z^\la\\=
87846 z^{41}+87846
   z^{39}+131769 z^{38}+87846
   z^{37}+146410 z^{36}+29282
   z^{35}+87846 z^{34}\\
+87846
   z^{33}+131769 z^{32}+123783
   z^{30}+146410 z^{29}+65219
   z^{28}+151734 z^{27}\\
+153065
   z^{26}+105149 z^{25}+154396
   z^{24}+145079 z^{23}+153065
   z^{22}+22627 z^{21}\\
+103818
   z^{20}+4719 z^{19}+78529
   z^{18}+156453 z^{17}+153186
   z^{16}+64614 z^{15}\\
+123178
   z^{14}+20933 z^{13}+154033
   z^{12}+84579 z^{11}+93533
   z^{10}+151492 z^9\\
+28325
   z^8+136730 z^7+23727
   z^6+43164 z^5+75636
   z^4+149358 z^3+126445
   z^2\\+97383
   z+7+\frac{80547}{1-z}+\frac{6
   809}{(1-z)^2}+\frac{17787}{
   (1-z)^3}+\frac{41261}{(1-z)
   ^4}+\frac{14641}{(1-z)^5}\\
\quad \text{\em modulo }11^5.
\end{multline}
\end{theorem}}

The next prime, $13$, is congruent to $1$ modulo~$6$, and
Theorem~\ref{thm:Freep} says that we may present the solution modulo
$13^e$ as a rational function with denominator being 
the $e$-th power of a quadratic factor.
It turns out that this quadratic factor is $(1-2z)(1+5z)$.
Experimentally, it seems that the maximal exponent $s$
such that $(1+5z)^{-s}$ appears in the (cleared) denominator
is actually $s=(e+1)/2$ (instead of $e$;
see for example Theorem~\ref{thm:Free13^5}), but our methods
are not able to demonstrate this.
What we are able to prove is summarised in the theorem below.

\begin{theorem} \label{thm:Free13}
Let $\al$ be a positive integer.
Then the generating function $F(z)=1+\sum_{\la\ge1}f_\la\,z^\la$, when 
coefficients are reduced modulo~$13^\al$, equals
$P_\al(z)/\big((1-2z)(1+5z)\big)^\al$,
where $P_\al(z)$ is a polynomial in $z$ over the integers. Moreover,
the sequence $(f_\la)_{\la\ge1}$ is ultimately periodic, with 
minimal period length equal to $12\cdot 13^{\al-1}$.  
\end{theorem}

The algorithm which led to the above theorem generates the following result 
for the modulus $13^5=371293$.

{
\allowdisplaybreaks
\begin{theorem} \label{thm:Free13^5}
We have
\begin{multline} \label{eq:free13^5}
1+\sum_{\la\ge1}f_\la\,z^\la\\=
314171 z^{42}+285610
   z^{40}+142805 z^{38}+114244
   z^{37}+285610 z^{36}+57122
   z^{35}
\kern1cm\\
+118638 z^{34}+285610
   z^{33}+325156 z^{32}+142805
   z^{30}+90077 z^{29}+338338
   z^{28}\\
+349323 z^{27}+188942
   z^{26}+103259 z^{25}+26364
   z^{24}+35152 z^{23}+188942
   z^{22}\\
+4732 z^{21}+76895
   z^{20}+310622 z^{19}+28561
   z^{18}+340535 z^{17}+358787
   z^{16}\\
+353379 z^{15}+135031
   z^{14}+115596 z^{13}+20280
   z^{12}+328874 z^{11}+55939
   z^{10}\\+116441 z^9+56745
   z^8+179309 z^7+342212
   z^6+219700 z^5+24336
   z^4\\
+238953 z^3+332462
   z^2+354965
   z+13\\+\frac{208033}{1+5
   z}+\frac{363181}{(1+5
   z)^2}+\frac{171366}{(1+5
   z)^3}
+\frac{334822}{1-2
   z}\\
+\frac{176228}{(1-2
   z)^2}+\frac{154635}{(1-2
   z)^3}+\frac{134017}{(1-2
   z)^4}+\frac{314171}{(1-2
   z)^5}
\quad \text{\em modulo }13^5.
\end{multline}
\end{theorem}}

For $p=17$,
the denominators of rational solutions to \eqref{eq:fmdiff} with $m=1$
(viewed modulo a power of $p$)
are powers of $1 + 15 z + 7 z^2$ (which does not 
factor modulo~$17$), and for $p=19,23$ they
are powers of (irreducible) ternary factors, etc. 
A generalisation of the argument in the proof of Theorem~\ref{thm:Free7}
shows that it is still true
that the period of $f_\la$ modulo a prime power $p^\al$ is a multiple 
of $p^{\al-1}$, but, which multiple it actually is, is in general
difficult to describe. For example, modulo~17 the period of $f_\la$
is $6\cdot 16$, modulo~$17^2$ it is $18\cdot 16\cdot 17$, 
while modulo~$17^3$ it is $102\cdot 16\cdot 17^2$.

\section{Free subgroup numbers for lifts of $\mathfrak{H}(4)$ modulo
prime powers}
\label{Sec:Freem4}

\noindent
In this last section, we apply Theorem~\ref{thm:ABCD} to the free 
subgroup numbers of another family of groups: the lifts $\Gamma_m(4)$
of the Hecke group $\mathfrak{H}(4)=C_2*C_4$ defined by
$$
\Gamma_m(4) = C_{2m} \underset{C_m}{\ast} C_{4m} = \big\langle
x,y\,\big\vert\, x^{2m} = y^{4m} = 1,\, x^2 = y^4\big\rangle,
$$
where $m$ is a positive integer. Given
$\lambda\geq1$, let $f_\la^{(4)}(m)$ denote the
number of free subgroups of index $4m\la$ in $\Gamma_m(4)$, and let
$F_m^{(4)}(z)=1+\sum_{\la\ge1}f_\la^{(4)}(m)z^\la$ 
be the corresponding generating function.
In particular, $f_\la^{(4)}:=f_\la^{(4)}(1)$ is the number of free 
subgroups of index $4\la$ in the Hecke group $\mathfrak{H}(4)$.

From the theory developed in \cite{MuHecke},
it is well-known that $F_m^{(4)}(z)$ satisfies a differential
equation. In order to determine the precise form of this differential
equation, several invariants of the groups $\Gamma_m(4)$ have to be
computed. Following the notation of \cite{MuHecke}, these are:
\begin{align*} \label{}
m_{\Gamma_m(4)} &= 4m,\\
\chi(\Gamma_m(4)) &= \frac{1}{2m} + \frac{1}{4m} - \frac{1}{m} 
    = -\frac{1}{4m},\\
\mu(\Gamma_m(4)) &= 1- m_{\Gamma_m(4)} \chi(\Gamma_m(4)) = 2,
\end{align*}
and the invariants $A_\mu = A_\mu(\Gamma_m(4))$, $0\leq \mu\leq 2$, 
for which we obtain
$$
A_0 = 3 m^2,\quad 
A_1 = 32 m^2,\quad 
A_2 = 16m^2.
$$
The differential equation then becomes\footnote{Cf.\ 
Formula~(18) in \cite[Sec.~2.3]{MuHecke} with $F_m^{(4)}(z) =
  1+z\beta_{\Gamma_m(4)}(z)$ for the result in a more general
  context.} 
\begin{equation} \label{eq:fmdiff4}
(1-(4m-2)z)F_m^{(4)}(z)-4mz^2(F_m^{(4)})'(z)-z(F_m^{(4)})^2(z)
-1-(1-4m+3m^2)z=0.
\end{equation}
Visibly, this is the special case of \eqref{eq:diffeqABCD} where
$A=4m-2$, $B=4m$, $C=1$, and  $D=1-4m+3m^2$. Hence,
Theorem~\ref{thm:ABCD} applies and allows us to derive results
analogous to the ones in Sections~\ref{Sec:Freem} and
\ref{Sec:Periods} for the lifts of the modular group. 
Since also the proofs are completely analogous, we content ourselves
with just stating these results.

\begin{theorem} \label{thm:Freem4}
Let $m$ be a positive integer,
$p$ a prime number with $p\ge 3,$ and let $\al$ be a positive integer. 
Then the generating function $F_m^{(4)}(z)=1+\sum_{\la\ge1}f_\la^{(4)}(m)z^\la,$
when coefficients are reduced modulo $p^\al,$ can be represented as a
rational function. Equivalently, the sequence $(f_\la^{(4)}(m))_{\la\ge1}$
is ultimately periodic modulo $p^\alpha$.
\end{theorem}

\begin{lemma} \label{lem4:Qnp1}
Let $p$ be a prime with $p\equiv1$~{\em(mod $4$)}.
For $A=4m-2,$ $B=4m,$ $C=1,$ $D=1-4m+3m^2$ {\em(}and, hence, $E=2m${\em),}
and $n\equiv\frac {p-1} {4},\frac {3(p-1)} {4}$~{\em(mod~$p$),} we have
$$
Q_n(z)=Q_{(p-1)/4}(z)\quad \text{\em modulo $p$},
$$
with $Q_n(z)$ as given in the assertion of Theorem~{\em\ref{thm:ABCD}}.
Furthermore, if $p\nmid m$, then 
the polynomial $Q_{(p-1)/4}(z)$ has degree $(p-1)/4$ in~$z$
when coefficients are reduced modulo~$p$, whereas 
$Q_{(p-1)/4}(z)=1$~{\em modulo}~$p$ if $p\mid m$.
\end{lemma}

\begin{lemma} \label{lem4:Qnp2}
Let $p$ be a prime with $p\equiv3$~{\em(mod $4$)}.
For $A=4m-2,$ $B=4m,$ $C=1,$ $D=1-4m+3m^2$ {\em(}and, hence, $E=2m${\em),}
and $n\equiv\frac {p-3} {4},\frac {3p-1} {4}$~{\em(mod~$p$),} we have
$$
Q_n(z)=Q_{(p-3)/4}(z)\quad \text{\em modulo $p$},
$$
with $Q_n(z)$ as given in Theorem~{\em\ref{thm:ABCD}}.
Furthermore, if $p\nmid m$, then 
the polynomial\break $Q_{(p-3)/4}(z)$ has degree $(p-3)/4$
in~$z$ when coefficients are reduced modulo~$p$, whereas 
$Q_{(p-3)/4}(z)=1$~{\em modulo}~$p$ if $p\mid m$.
\end{lemma}

\begin{theorem} \label{thm:Freep4}
Let $\al$ be a positive integer.
If $p$ is a prime with $p\equiv1$~{\em(mod $4$),} 
then the generating function 
$F_m^{(4)}(z)=1+\sum_{\la\ge1}f_\la^{(4)}(m)\,z^\la,$ when 
coefficients are reduced modulo~$p^\al,$ equals 
$\overline P_\al(z)/Q_{(p-1)/4}^\al(z),$
where $\overline P_\al(z)$ is a polynomial in $z$ over the integers,
and $Q_{(p-1)/4}(z)$ is the polynomial of degree $(p-1)/4$
given by the formula in Theorem~{\em\ref{thm:ABCD}} with
$A=4m-2,$ $B=4m,$ $C=1,$ $D=1-4m+3m^2$ and $E=2m$.

On the other hand,
if $p$ is a prime with $p\equiv3$~{\em(mod $4$),} 
then the generating\break function 
$F_m^{(4)}(z)=1+\sum_{\la\ge1}f_\la^{(4)}(m)\,z^\la,$ when 
coefficients are reduced modulo~$p^\al,$ equals\break
$\widehat P_\al(z)/Q_{(p-3)/4}^\al(z),$
where $\widehat P_\al(z)$ is a polynomial in $z$ over the integers,
and $Q_{(p-3)/4}(z)$ is the polynomial of degree $(p-3)/4$
given by the formula in Theorem~{\em\ref{thm:ABCD}} with
$A=4m-2,$ $B=4m,$ $C=1,$ $D=1-4m+3m^2$ and $E=2m$.
\end{theorem}

\begin{remarks}
(1) Also in the present context, it would be possible without any
additional effort to make 
Theorem~\ref{thm:Freep4} more explicit for concrete primes~$p$; thus
obtaining results analogous to the ones in
Theorems~\ref{thm:Free7}--\ref{thm:Free13^5}. We
omit this here for the sake of brevity.

\medskip
(2) The Riccati differential equation \eqref{eq:fmdiff4} is in a form
so that the method from \cite{KKM} for the determination of the
behaviour of sequences modulo powers of $2$ as well as the method
from \cite{KrMu3Power} for the determination of the
behaviour of sequences modulo powers of $3$ are applicable.
Consequently, one could derive results analogous to the ones
in \cite[Sec.~8]{KKM} and \cite[Sec.~16]{KrMu3Power} for the free
subgroup numbers of the lifts $\Gamma_m(3)$ of $\PSL_2(\Z)$.

\medskip
(3) One may speculate that the obvious generalisation of
Theorems~\ref{thm:Freem} and \ref{thm:Freem4} may hold for
the lifts $\Gamma_m(q)$ of an arbitrary Hecke group
$\mathfrak{H}(q)$, given by
$$
\Gamma_m(q) = C_{2m} \underset{C_m}{\ast} C_{qm} = \big\langle
x,y\,\big\vert\, x^{2m} = y^{qm} = 1,\, x^2 = y^q\big\rangle,
\quad m\ge1. 
$$
Computer experiments seem to support this speculation.
\end{remarks}

\appendix
\section*{Appendix: How were the expressions for
$p_{n,k}$ and $q_{n,k}$ found?}

\setcounter{equation}{0}%
\global\def\theequation{\mbox{A.\arabic{equation}}}

When papers get written, the path(s) towards the results presented 
often get lost, are covered up, or --- in some cases
--- even hidden intentionally. As we believe that, in the present
case, this path is actually 
quite interesting and instructive, 
we shall describe in this appendix how we were led to conjecture the
formulae in \eqref{eq:pnk} and \eqref{eq:qnk}, without which we
would not even have been able to prove the mere existence of
solutions to \eqref{eq:idABCD}, and without which we would 
therefore never have reached the result described in Theorem~\ref{thm:Freem}
concerning the number of free subgroups in lifts of the inhomogeneous modular
group, which was our primary goal and motivation.

Recall  that the generating function 
$F(z)=1+\sum_{\la\ge1}f_\la\,z^\la$, where $f_\la$ denotes again
the number of free subgroups of $\PSL_2(\Z)$
of index $6\la$, satisfies the 
Riccati differential equation \eqref{eq:fmdiff} with $m=1$, that is, 
\begin{equation} \label{eq:FreeEF}
(1-4z)F(z)-6z^2F'(z)-zF^2(z)-1=0,
\end{equation}
This equation allows us to compute the numbers $f_\la$ recursively.
We started out by looking at tables of these numbers,
when they are reduced modulo powers of $5$, $7$, $11$, $13$, etc.
The case of $5$ is a trivial (and exceptional) one, since it
is not hard to see that $f_\la$ vanishes modulo any power of
$5$ once $\la$ is large enough. The next case, powers of $7$, is more
interesting. Computer experiments quickly led to the conjecture that
the free subgroup numbers $f_\la$ modulo a given power of $7$ are ultimately
periodic, and it was even possible to predict the period length.
Theorem~\ref{thm:Free7} presents the precise corresponding
statement. Continuing with powers of $11$, a similar picture
emerged. The precise statement is presented in
Theorem~\ref{thm:Free11}. With some effort, we could see that
the free subgroup numbers $f_\la$ are also periodic when reduced
modulo~$13$ or $17$, however with considerably larger period
lengths. What happens modulo higher powers of $13$ or $17$ was
unclear at that point.

We then decided to try to prove, say, the statements modulo powers of
$7$ and $11$.
An approach in the style of \cite{KKM} and \cite{KrMu3Power} came quite
close, but did not succeed. At some point, we had the idea of
``Pad\'e-approximating" the solution to \eqref{eq:FreeEF} by a rational
function with numerator and denominator degrees equal to $n$, and see how far
this rational function fails to satisfy the differential equation 
\eqref{eq:FreeEF}.
We were (pleasantly) surprised to see that the error on the
right-hand side was 
a ``round"\footnote{According to Hardy \cite[Ch.~III]{HardAA},
``a number is described in popular language as {\it round\/}
if it is the product of a considerable number of comparatively
small factors," which was picked up by Kuperberg \cite[Sec.~VII.A]{KupeAG}
(in fact, being unaware of Hardy's statement) in the form
``an integer is {\it round\/} if it is a product of relatively small
(prime) numbers. A round enumeration in
combinatorics almost always has an explicit product formula."
See \cite[Sec.~2]{KratBZ} for a more detailed elaboration of the
``notion" of ``roundness."} 
multiple of $z^{2n+1}/Q_n^2(z)$, 
and, using the
guessing programme {\tt Rate},\footnote{\label{foot:Rate}written in 
{\sl Mathematica} 
by the first author; available from\newline {\tt
http://www.mat.univie.ac.at/\lower0.5ex\hbox{\~{}}kratt}}
it was not difficult to come up with the conjecture that
\eqref{eq:Rmdiff} (with $m=1$) should hold.

After several failed attempts at proving this conjecture, in a 
mood of despair, we decided to make everything more general, not expecting
anything but horrendous results. We put as many parameters as possible
into the differential equation \eqref{eq:FreeEF}, so that we were led 
to consider the differential equation \eqref{eq:diffeqABCD}. 
We were in for another pleasant surprise: although one cannot go very far
with nowadays' computers (we only made it up to $n=5$), the pattern
which is displayed in what is now Theorem~\ref{thm:ABCD} became
quickly clear, however without any clue about the explicit
form of the rational function $R_n(z)$. 
For example, we obtained
\begin{gather} \label{eq:R1}
R_1(z)=\frac{1+
   (-B-C+D)z}
   {1+(-A-B-2
   C)z},\\
\label{eq:R2}
R_2(z)=\frac{1+(-A-4
   B-3
   C+D)z+\left(-A
   D+2 B^2+3 B
   C-3 B
   D+C^2-3
   C
   D\right)z^2}{1+ (-2 A-4
   B-4
   C)z+
   \left(A^2+3 A B+3 A
   C+2 B^2+6 B
   C+3
   C^2-C
   D\right)z^2},\\
\label{eq:R3}
R_3(z)=\frac {P_3(z)} {Q_3(z)},
\end{gather}
with
\begin{multline*}
P_3(z)=1
+z (-2
   A-9 B-5
   C+D)\\
+z^2 \left(A^2+7 A B+4 A
   C-2 A D+18
   B^2+22 B C-8 B
   D+6 C^2-6
   C
   D\right)\\
+z^3
   \big(A^2 D+6 A B
   D+4 A C
   D-6 B^3-11 B^2
   C+11 B^2
   D\\
-6 B
   C^2+18 B C
   D-C^3+6
   C^2
   D-C
   D^2\big)
\end{multline*}
and
\begin{multline*}
Q_3(z)=1
+z (-3 A-9
   B-6 C)\\
+z^2
   \left(3 A^2+15 A B+10 A
   C+18 B^2+30 B
   C+10 C^2-2
   C
   D\right)\\
+z^3
   \big(-A^3-6 A^2 B-4 A^2
   C-11 A B^2-18 A B
   C-6 A C^2+2
   A C D-6
   B^3\\
-22 B^2 C-18 B
   C^2+6 B C
   D-4 C^3+4
   C^2
   D\big).
\end{multline*}
On the basis of these computer data (including the ones for $R_4(z)$
and $R_5(z)$), it is not difficult to come up with guesses for the
first few coefficients $p_{n,k}$ and $q_{n,k}$:
{\allowdisplaybreaks
\begin{align}
\label{eq:pn1}
p_{n,1}&=-(n-1)A-n^2B-(2 n-1)C+D,\\
\label{eq:qn1}
q_{n,1}&=-nA-n^2B-2 nC,\\
\notag
p_{n,2}&=\tfrac{1}{2}  (n-2)
   (n-1)A^2
+\tfrac{1}{2}  (n-2)
   (n-1) (2 n+1)A B
+2    (n-2) (n-1)A C\\
\notag
&\kern1cm
-   (n-1)A D
+\tfrac{1}{2} 
   (n-1)^2 n^2B^2
+
   (n-1) \left(2 n^2-2
   n-1\right)B C\\
\label{eq:pn2}
&\kern1cm
-
   (n-1) (n+1)B D
+
   (n-1) (2 n-3)C^2
-3  (n-1)C
   D,\\
\notag
q_{n,2}&=\tfrac{1}{2}  (n-1)
   nA^2
+\tfrac{1}{2}  (n-1) n
   (2 n-1)A B
+ (n-1)
   (2 n-1)A C
+\tfrac{1}{2} 
   (n-1)^2 n^2B^2\\
\label{eq:qn2}
&\kern1cm
+
   (n-1) n (2 n-1)B C
+
   (n-1) (2 n-1)C^2
- (n-1)C
   D,\\
\notag
p_{n,3}&=-\tfrac{1}{6}  (n-3) (n-2)
   (n-1)A^3
-\tfrac{1}{2} 
   (n-3) (n-2)
   \left(n^2-n-1\right)A^2 B\\
\notag
&\kern1cm
-\tfrac{
   1}{2}  (n-3) (n-2) (2
   n-3)A^2 C
+\tfrac{1}{2} 
   (n-2) (n-1)A^2 D\\
\notag
&\kern1cm
-\tfrac{1}{6}  (n-3) (n-2) \left(3
   n^3-3
   n^2-n-2\right)A
   B^2\\
\notag
&\kern1cm
-\tfrac{1}{2}
    (n-3) (n-2) (2 n-3)
   (2 n+1)A B C
+\tfrac{1}{2} 
   (n-2) (n+1) (2 n-3)A B D\\
\notag
&\kern1cm
-
   (n-3) (n-2) (2 n-3)A C^2
+
   (n-2) (3 n-5)A C D\\
\notag
&\kern1cm
-\tfrac{1}{6}
    (n-2)^2 (n-1)^2
   n^2B^3
-\tfrac{1}{6} (n-2)
   (2 n-3) \left(3 n^3-6
   n^2-n-2\right) B^2 C\\
\notag&\kern1cm
+\tfrac{1}{2}
   (n-2)
   \left(n^3-n-2\right) B^2 D
-
   (n-2) (2 n-3) \left(n^2-2
   n-1\right)B C^2\\
\notag&\kern1cm
+ (n-2)
   \left(3 n^2-2
   n-3\right)B C D
-\tfrac{1}{3} 
   (n-2) (2 n-5) (2 n-3)C^3\\
\label{eq:pn3}
&\kern1cm
+2 
   (n-2) (2 n-3)D C^2
- (n-2)C D^2,\\
\notag
q_{n,3}&=-\tfrac{1}{6}  (n-2) (n-1)
   nA^3
-\tfrac{1}{2}  (n-2)
   (n-1)^2 nA^2 B
- (n-2)
   (n-1)^2A^2 C\\
\notag
&\kern1cm
-\tfrac{1}{6} 
   (n-2) (n-1) n \left(3 n^2-6
   n+2\right)A B^2
- (n-2)
   (n-1) n (2 n-3)A B C\\
\notag
&\kern1cm
- (n-2)
   (n-1) (2 n-3)A C^2
+ (n-2)
   (n-1)A C D\\
\notag
&\kern1cm
-\tfrac{1}{6} 
   (n-2)^2 (n-1)^2
   n^2B^3
-\tfrac{1}{3}  (n-2)
   (n-1) n \left(3 n^2-6
   n+2\right)B^2 C\\
\notag
&\kern1cm
- (n-2)
   (n-1) n (2 n-3)B C^2
+ (n-2)
   (n-1) nB C D\\
\label{eq:qn3}
&\kern1cm
-\tfrac{2}{3} 
   (n-2) (n-1) (2 n-3)C^3
+2 
   (n-2) (n-1)C^2 D.
\end{align}}%
So, it seemed ``clear" that the $p_{n,k}$'s and $q_{n,k}$'s were homogeneous 
polynomials in $A,B,C,\break D$ of degree $k$.
However, having stared at this for some while, it seemed hopeless to us 
to continue along these lines.

We made the ``intermediate" observation that, should we be able to
demonstrate something at all in this context, it would have to come from the
``hypergeometric world," and consequently the quadratic factors in the
right-hand side product of \eqref{eq:idABCD} must factor into linear
factors. This happens only if the discriminant
of the factor, viewed as a quadratic form in $\ell$, has a square root, and
the latter requires the parametrisation $E^2=A^2-4CD$. If one uses
this parametrisation, then the right-hand side product in
\eqref{eq:idABCD} becomes
$$
- (A+C+D)
\prod _{\ell=1} ^{n}(\ell AB+AC+CD+\ell^2B^2+2\ell BC+C^2)
=-\frac {1} {C}\Pi_+\Pi_-,
$$
where, again, we make use of the short-hand notation in \eqref{eq:Pi}.

At a certain point, we decided to ``look at the other end," that is,
to try to find formulae for the coefficients $p_{n,k}$ and $q_{n,k}$
for $k=n,n-1,\dots$ Comparing coefficients of $z^{2n+1}$ on both
sides of \eqref{eq:idPQ}, one obtains the equation
$$
-Ap_{n,n}q_{n,n}
-Cp^2_{n,n}
-Dq^2_{n,n}=
-\frac {1} {C}\Pi_+\Pi_-.
$$
Moreover, upon using the parametrisation $E^2=A^2-4CD$, the left-hand
side factors, leading to
$$
\left(Cp_{n,n}+\tfrac {1} {2}(A-E)q_{n,n}\right)
\left(Cp_{n,n}+\tfrac {1} {2}(A+E)q_{n,n}\right)=
\Pi_+\Pi_-.
$$
If there is any ``justice" in the world, then one of the factors on
the left-hand side 
would have to agree with one of the factors on the right-hand side, and the
same should hold for the other two factors. Indeed, the computer told us
that, apparently,
\begin{align*}
Cp_{n,n}+\tfrac {1} {2}(A-E)q_{n,n}&=\Pi_-,\\
Cp_{n,n}+\tfrac {1} {2}(A+E)q_{n,n}&=\Pi_+.
\end{align*}
Solving this system of linear equations 
for $p_{n,n}$ and $q_{n,n}$, one arrives at
$$
p_{n,n}=\frac {(-1)^{n+1}} {2CE}
\big(    (A  -  E)\Pi_+
-(A+E)\Pi_-
\big)
$$
(that is, Equation~\eqref{eq:pnn}), and 
$$
q_{n,n}=\frac {(-1)^{n}} {E}
(    \Pi_+-\Pi_-)
$$
(that is, Equation~\eqref{eq:qnn}).

Subsequently, we proceeded to compare coefficients of $z^{2n}$ on both
sides of \eqref{eq:idPQ}. This leads to a linear equation
for $p_{n,n-1}$ and $q_{n,n-1}$, with coefficients containing $\Pi_+$
and $\Pi_-$. Containing two unknowns, it is certainly undetermined.
Nevertheless, one is able to come up with a guess:
one considers the equation first modulo $\Pi_+-\Pi_-$, and then
modulo $\Pi_++\Pi_-$. As it turns out, this leads to dramatic
simplification, and one obtains two {\it linear} congruences for
$p_{n,n-1}$ and $q_{n,n-1}$ (the moduli being $\Pi_+-\Pi_-$ and
$\Pi_++\Pi_-$, respectively), with coefficients being rational
functions in $A,B,C,E$. One can then again use computer experiments to
come up with a guess for the multiplicative factors which make these
congruences into equations. Once this is done, one has two linear
equations, which one solves for $p_{n,n-1}$ and $q_{n,n-1}$. The
result was \eqref{eq:pnk} and \eqref{eq:qnk} with $k=1$ (albeit not
yet written in this form). We did the same for $p_{n,n-2}$ and
$q_{n,n-2}$, and for $p_{n,n-3}$ and $q_{n,n-3}$. (The latter required
to also compute $R_6(z)$ and $R_7(z)$, which now, by putting the
knowledge of the coefficients $p_{n,k}$ and $q_{n,k}$ for
$k=n,n-1,n-2$ into the computer programme, became feasible.)

At this stage, we realised that the coefficients of $\Pi_+$ and
$\Pi_-$ became ``elegant," if one expanded them --- viewed as polynomials
in $n$ --- in the basis
$$
1,\ (n-k+1),\ (n-k+1)(n-k+2),\ \dots,\ 
(n-k+1)(n-k+2)\cdots (n-1)n.
$$
We also had at this stage abundant ``information" on the arising
denominators and signs.
More precisely, our experimental data told us that one could
apparently write
\begin{align*} 
p_{n,n-k}&=\frac {(-1)^{n+1}B^{n-k}} {2C\,(\frac {E} {B}-k)_{2k+1}}\\
&\kern.5cm
\times
\Bigg( \left(\tfrac {A + 2 C + E} {2B}\right)_{n+1}
\sum_{j=0}^k (n-k+1)_{k-j}
\left(-\tfrac E B
  + j + 1\right)_{ 
    k - j} \left(\tfrac {A + 2 C - E} {2B}\right)_j\\
&\kern5cm
\cdot
    \big (p^{(1)}_{n,k,j} A +
p^{(2)}_{n,k,j} B - p^{(3)}_{n,k,j}E\big)\\
&\kern1.5cm 
- \left(\tfrac {A + 2 C - E} {2B}\right)_{n+1}
\sum_{j=0}^k (n-k+1)_{k-j}
\left(\tfrac E B
  + j + 1\right)_{ 
    k - j} \left(\tfrac {A + 2 C + E} {2B}\right)_j\\
&\kern5cm
\cdot
    \big (p^{(1)}_{n,k,j} A +
p^{(2)}_{n,k,j} B + p^{(3)}_{n,k,j}E\big)
\Bigg)
\end{align*}
and
\begin{align} 
\notag
q_{n,n-k}&=\frac {(-1)^{n}B^{n-k}} {(\frac {E} {B}-k)_{2k+1}}\\
\notag
&\kern.5cm
\times
\Bigg( \left(\tfrac {A + 2 C + E} {2B}\right)_{n+1}
\sum_{j=0}^k q_{n,k,j} (n-k+1)_{k-j}
\left(-\tfrac E B
  + j + 1\right)_{ 
    k - j} \left(\tfrac {A + 2 C - E} {2B}\right)_j\\
\label{eq:qAnsatz}
&\kern1.5cm 
- \left(\tfrac {A + 2 C - E} {2B}\right)_{n+1}
\sum_{j=0}^k q_{n,k,j} (n-k+1)_{k-j}
\left(\tfrac E B
  + j + 1\right)_{ 
    k - j} \left(\tfrac {A + 2 C + E} {2B}\right)_j
\Bigg),
\end{align}
with ``nice" numerical coefficients $p^{(1)}_{n,k,j},p^{(2)}_{n,k,j},
p^{(3)}_{n,k,j},q_{n,k,j}$. The available data at this point (namely,
the rational functions $R_n(z)$ for $n=1,2,\dots,7$) were largely 
insufficient to make any guesses about these coefficients. However,
assuming the above forms of $p_{n,k}$ and $q_{n,k}$, together with
the earlier observed polynomiality of $p_{n,k}$ and $q_{n,k}$
in $A,B,C,D$, and, hence, also in $A,B,C,E$
(except for possible powers of $C$ in denominators), sufficient data could 
now be easily produced! Let us first concentrate on $q_{n,n-k}$ (as 
we did). If $q_{n,n-k}$ is to be a polynomial, then
$(\frac {E} {B}-k)_{2k+1}$ must divide the expression between parentheses
on the right-hand side of \eqref{eq:qAnsatz}, that is, this expression
must vanish for $E=Bs$, $s=-k,-k+1,\dots,k$. Let us do the
substitution $E=Bs$ in this expression, and let us suppose that $s>0$.
Then we have
$$ \left(\tfrac {A + 2 C + E} {2B}\right)_{n+1}=
 \left(\tfrac {A + 2 C + Bs} {2B}\right)_{n+1}
=\left(\tfrac {A + 2 C + Bs} {2B}\right)_{n+1-s}
\left(\tfrac {A + 2 C - Bs} {2B}+n+1\right)_{s}.
$$
Comparing with
$$ \left(\tfrac {A + 2 C - E} {2B}\right)_{n+1}=
 \left(\tfrac {A + 2 C - Bs} {2B}\right)_{n+1},
$$
we infer that $\left(\tfrac {A + 2 C - Bs} {2B}+n+1\right)_{s}$ must
divide the second sum over $j$ {\it as a polynomial in $n$}.
This provides many non-trivial vanishing conditions for this sum over
$j$, viewed as polynomial in $n$ (they actually ``overdetermine" it), 
and this suffices to compute 
the coefficients $q_{n,k,j}$ for a large range of $k$'s and
corresponding $j$'s. Using the guessing programme {\tt Rate} (see
Footnote~\ref{foot:Rate}) once again, it is then straightforward to guess
the formula in \eqref{eq:qnk}. For the coefficients 
$p^{(1)}_{n,k,j},p^{(2)}_{n,k,j},
p^{(3)}_{n,k,j}$ one does not have to do the same work again:
it suffices to look at the quotients
$$
\frac {p^{(1)}_{n,k,j}} {q_{n,k,j}},\ 
\frac {p^{(2)}_{n,k,j}} {q_{n,k,j}},\ 
\frac {p^{(3)}_{n,k,j}} {q_{n,k,j}},
$$
which quickly leads one to the formula in \eqref{eq:pnk}.


\begin{thebibliography}{99}

\bibitem{BailAA}
W.\,N.~Bailey, {\em Generalized Hypergeometric Series},
Cambridge University Press, Cambridge, 1935.

\bibitem{BaGrAA}
G.\,A. Baker, Jr. and P.\,R. Graves-Morris, {\em Pad\'e Approximants}, 
Cambridge University Press, Cambridge, 1996.

\bibitem{BundAA}
P. Bundschuh, {\em Einf\"uhrung in die Zahlentheorie},
Springer-Verlag, Berlin, 1992.



\bibitem{CaMuDescent} P.\,J. Cameron and T.\,W. M\"uller,
A descent principle in modular subgroup arithmetic, {\it J. Pure
Appl.\ Algebra} {\bf 203} (2005), 189--203.

\bibitem{Frob1} G. Frobenius, Verallgemeinerung des Sylow'schen Satzes, 
{\it Sitz.ber.\ K\"onigl.\ Preuss.\ Akad.\ Wiss.\ Berlin}, 1895, 981--993.

\bibitem{Frob2} G. Frobenius, \"Uber einen Fundamentalsatz der
Gruppentheorie, {\it
  Sitz.ber. K\"onigl.\ Preuss.\ Akad.\ Wiss.\ Berlin}, 1903,
987--991. 

\bibitem{GospAB} 
R.\,W. Gosper, Decision procedure for 
indefinite hypergeometric summation, {\it Proc. Natl. Acad. Sci.\
USA} {\bf 75} (1978), 40--42.

\bibitem{GrKPAA} R.\,L.~Graham, D.\,E.~Knuth and O.~Patashnik,
{\it Concrete Mathematics}, Addison-Wesley, Reading,
Massachusetts, 1989.

\bibitem{PHall1} P. Hall, A contribution to the theory of groups of
  prime-power order, \textit{Proc. Lond. Math. Soc.} (2) \textbf{36}
  (1934), 29--95.  

\bibitem{PHall2} P. Hall, On a theorem of Frobenius,
  \textit{Proc.\ London Math.\ Soc.} (2) \textbf{40} (1935),
  468--501. 

\bibitem{HardAA}
G. H. Hardy, {\it Ramanujan}, twelve lectures on subjects suggested
by his life and work, Chelsea Publ.\ Comp., New York, 1940.

\bibitem{KKM} 
M. Kauers, C. Krattenthaler and T.\,W. M\"uller, A
method for determining the mod-$2^k$ behaviour of
recursive sequences, with applications to subgroup counting,
{\it Electron. J. Combin.} {\bf 18}(2) (2012), Article~P37, 83~pp.

\bibitem{KratBZ}
C.    Krattenthaler, Advanced determinant calculus: a complement,
{\it Linear Algebra Appl.} {\bf 411} (2005), 68--166.

\bibitem{KratMuHecke} 
C. Krattenthaler and T.\,W. M\"uller,
Parity patterns associated with lifts of Hecke groups,
{\it Abh.\ Math.\ Sem.\ Univ.\ Hamburg} {\bf 78} (2008), 99--147.

\bibitem{KrMu3Power} 
C. Krattenthaler and T.\,W. M\"uller, A
method for determining the mod-$3^k$ behaviour of
recursive sequences, preprint, 2012.

\bibitem{KupeAG}
G.    Kuperberg, An exploration of the permanent-determinant method,
{\it Electron. J. Combin.} {\bf 5} (1998), Article~\#R46, 34~pp.

\bibitem{LegeAA}
A.\,M. Legendre, {\em Essai sur la th\'eorie des nombres}, 
2$^{\text{e}}$ ed., 
Courcier, Paris, 1808.

\bibitem{MuHecke} T.\,W. M\"uller, Parity patterns in Hecke groups and
  Fermat primes. In: {\em Groups: Topological, Combinatorial, and
    Arithmetic Aspects,} Proceedings of a conference held 1999 in
  Bielefeld (T.\,W.~M\"uller, ed.), LMS Lecture Note Series, vol.~311,
  Cambridge University Press, 
Cambridge, 2004, pp.~327--374.

\bibitem{MuDescent} T.\,W. M\"uller, Modular subgroup arithmetic and a
theorem of Philip Hall,
{\it Bull.\ London Math.\ Soc.} {\bf 34} (2002), 587--598.

\bibitem{MuSubArith} T.\,W. M\"uller, Modular subgroup arithmetic in free
products, {\it Forum Math.} {\bf 15} (2003), 759--810.

\bibitem{MuPu2} T.\,W. M\"uller and J.-C. Schlage-Puchta, 
Modular arithmetic of free subgroups, 
{\it Forum Math.} {\bf 17} (2005), 375--405.

\bibitem{MuPu} T.\,W. M\"uller and J.-C. Schlage-Puchta, Divisibility
properties of subgroup numbers for the modular group, \textit{New York
J. Math.} \textbf{11} (2005), 205--224. 

\bibitem{Ore} O. Ore, \textit{Number Theory And Its History},
Dover Publ.\ Inc., New York, 1988.

\bibitem{PeWZAA} M.~Petkov\v sek, H.~Wilf and D.~Zeilberger,
{\em A=B}, A.~K.~Peters, Wellesley, 1996.

\bibitem{SlatAC} L.\,J.~Slater, {\em Generalized Hypergeometric Functions},
Cambridge University Press, Cambridge, 1966.

\bibitem{Sylow}
L. Sylow, Th\'eor\`emes sur les groupes de substitutions, 
{\it Math. Ann.} {\bf 5} (1872), 584--594.

\end{thebibliography}
\end{document}